\documentclass[12pt,amsfonts,amscd]{amsart}

\usepackage{amssymb,amsmath,amsthm,longtable,multirow,array}
\usepackage{graphicx}
\usepackage{datetime}
\usepackage{amsfonts}
\usepackage{color}
\usepackage{mathrsfs}
\newtheorem{Theorem}{Theorem}[section]
\newtheorem{Proposition}[Theorem]{Proposition}

\newtheorem{Lemma}[Theorem]{Lemma}
\newtheorem{Corollary}[Theorem]{Corollary}
\theoremstyle{definition}\newtheorem{Definition}[Theorem]{Definition}
\newtheorem{Example}[Theorem]{Example}
\newtheorem{Exercise}[Theorem]{Exercise}

\newtheorem{Para}[Theorem]{}
\theoremstyle{remark}
\newtheorem*{remark*}{Remark}

\providecommand\ba[1]{\begin{align*}#1\end{align*}}
\providecommand\baa[1]{\begin{align}#1\end{align}}
\providecommand\baaa[1]{\begin{equation}\begin{split}#1\end{split}\end{equation}}

\providecommand\brs{\begin{remark*}}
\providecommand\ers{\end{remark*}}

\providecommand\be{\begin{enumerate}}
\providecommand\ee{\end{enumerate}}
\providecommand\bT{\begin{Theorem}}
\providecommand\eT{\end{Theorem}}
\providecommand\bP{\begin{Proposition}}
\providecommand\eP{\end{Proposition}}
\providecommand\bD{\begin{Definition}}
\providecommand\eD{\end{Definition}}
\providecommand\bE{\begin{Example}}
\providecommand\eE{\end{Example}}
\providecommand\bEE{\begin{Exercise}}
\providecommand\eEE{\end{Exercise}}
\providecommand\bPP{\begin{Para}}
\providecommand\ePP{\end{Para}}
\providecommand\bL{\begin{Lemma}}
\providecommand\eL{\end{Lemma}}
\providecommand\bC{\begin{Corollary}}
\providecommand\eC{\end{Corollary}}
\providecommand\bpp{\begin{proof}}
\providecommand\epp{\end{proof}}
\providecommand\bee{\begin{equation}}
\providecommand\eee{\end{equation}}

\providecommand\beqq{\begin{eqnarray*}}
\providecommand\eeqq{\end{eqnarray*}}

\providecommand\bay{\begin{array}}
\providecommand\eay{\end{array}}

\providecommand\epf{\hskip.2in\vrule width.4pt height6pt depth0pt\vrule
width5.2pt height6pt depth-5.6pt\hskip-5.2pt\vrule width5.2pt
height.4pt depth0pt\vrule width.4pt height6pt depth0pt\ }

\providecommand\rk{\mathcal}

\providecommand\CC{{\Bbb C}}
\providecommand\RR{{\Bbb R}}
\providecommand\ZZ{{\Bbb Z}}

\providecommand\NN{{\Bbb N}}

\providecommand\ol{\overline}

\providecommand\al{\alpha}

\providecommand\dl{\delta}

\providecommand\eps{\varepsilon}
\providecommand\Gm{\Gamma}
\providecommand\gm{\gamma}

\providecommand\Lm{\Lambda}
\providecommand\lm{\lambda}

\providecommand\Om{\Omega}
\providecommand\om{\omega}

\providecommand\sbs{\subset}

\providecommand\arg{\operatorname{arg}}

\providecommand\Id{\operatorname{Id}}

\providecommand\diam{\operatorname{diam}}

\providecommand\dis{\operatorname{dis}}
\providecommand\Cr{\operatorname{Cr}}

\providecommand\iiff{if and only if }

\begin{document}

\title[Strict Whitney arcs and $t$-quasi  self-similar arcs]{On strict Whitney arcs and $t$-quasi  self-similar arcs}
\author{Daowei Ma}
\author{Xin Wei}
\author{Zhi-Ying Wen}

\begin{abstract} A connected compact subset $E$ of $\RR^N$ is said to be a strict Whitney set if there exists a real-valued $C^1$ function $f$ on $\RR^N$ with $\nabla f|_E\equiv 0$ such that $f$ is constant on no non-empty relatively open subsets of $E$. We prove that each  self-similar arc of Hausdorff dimension $s>1$ in $\RR^N$ is a strict Whitney set with criticality $s$. We also study a special kind of self-similar arcs, which we call ``regular'' self-similar arcs. We obtain necessary and sufficient conditions for a regular self-similar arc $\Lm$ to be a $t$-quasi-arc, and for  the Hausdorff measure function on $\Lm$ to be a strict Whitney function. We prove that if a regular self-similar arc has ``minimal corner angle'' $\theta_{\min}>0$, then it is a 1-quasi-arc and hence its Hausdorff measure function is a strict Whitney function. We provide an example of a one-parameter family of regular self-similar arcs with various features. For some value of the parameter $\tau$, the Hausdorff measure function of the self-similar arc is a strict Whitney function on the arc, and hence the self-similar arc is an $s$-quasi-arc, where $s$ is the Hausdorff dimension of the arc.  For each $t_0\ge 1$, there is a value of $\tau$ such that the corresponding self-similar arc is a $t$-quasi-arc for each $t>t_0$, but it is not a $t_0$-quasi-arc. For each $t_0>1$, there is a value of $\tau$ such that the corresponding self-similar arc is a $t_0$-quasi-arc, but it is a $t$-quasi-arc for no $t\in [1, t_0)$.
\end{abstract}

\address{dma@math.wichita.edu, Department of Mathematics, Wichita
State University, Wichita, KS 67260-0033, USA}
\address{wei@math.wichita.edu, Department of Mathematics, Wichita
State University, Wichita, KS 67260-0033, USA}
\address{wenzy@mail.tsinghua.edu.cn, Department of Mathematics, Tsinghua University, Beijing, 100080, P.R.China}

\maketitle

\bigskip

\section{Introduction}

In fractal geometry, Morse-Sard Theorem (see \cite{Sa}) states that if $f \in C^k(\RR^m,\RR^N)$ with $k \ge\max(m-N+1,1)$, then the set of critical values of $f$ has zero Lebesgue measure in $\RR^N$. However, Whitney in 1935 constructed a differentiable function $f:\RR^2\to\RR$ whose critical set is a fractal planar arc $\gm$ with Hausdorff dimension $\log4/\log3$, and whose set $f(\gm)$ of critical values contains an interval and therefore has positive Lebesgue measure (see \cite{Whi2}). This is called Whitney phenomenon; it seems to contradict the Morse-Sard Theorem. It is due to the fact that the arc $\gm$ is a fractal and $f$ has lower smoothness. Such a set is called a Whitney set.
  
\bD \label{df1.1}
{\rm A connected set $E \subset \RR^N$ is said to be a {\it Whitney set}, if there is a $C^1$ function $f:\RR^N \to \RR$ such that $\nabla f|_E\equiv 0$ but $f|_E$ is not constant. The function $f$ is said to be a {\it Whitney function} for $E$, and its restriction $f|_E$ to $E$ is said to be a {\it Whitney function} on $E$. If a Whitney function $f|_E$ on $E$ is non-constant on each non-empty relatively open subset of $E$, then $f|_E$ is said to be  a {\it strict Whitney function} on $E$, $f$ is said to be a {\it strict Whitney function} for $E$,  and the set  $E$ is said to be  a {\it strict Whitney set}.} 
\eD

The following special case of the Whitney Extension Theorem \cite{Whi} will be used.

\bL \label{lm3.1}
Suppose that $E\subset \RR^N$ is compact and $f: E \to \RR$ is a function. If for each $\eps >0$, there exists $\dl >0$ such that for each pair of points $x,y \in E$ with $|x-y|<\dl$, one has $|f(x)-f(y)|\le\eps|x-y|$, then there is a $C^1$ extension $\tilde{f}:\RR^N \to \RR$ of $f$ such that $\tilde{f}|_E=f$ and $\nabla \tilde{f}|_E\equiv 0$.
\eL

Lemma~\ref{lm3.1} suggests the following definition.

\bD \label{abs}{\rm A compact connected metric space $A$ is said to be a {\it Whitney set} if there is a non-constant function $f: A\to\RR$ such that $|f(x)-f(y)|=o(d_A(x,y))$ for $x,y\in A$.}
\eD

By Lemma~\ref{lm3.1}, for a compact connected subset $A$ of $\RR^N$, Definition~\ref{abs} is consistent with Definition~\ref{df1.1}.

About Whitney sets, we know the following.

(a) For a set $E\subset \RR^N$, if every pair of points in $E$ are connected by a rectifiable arc lying in $E$, then $E$ is not a Whitney set (Whyburn \cite{Wh}, 1929).

(b) For a continuous function $g:\RR \to \RR$, the graph $G$ of $g$ is not a Whitney set (Choquet \cite{Ch}, 1944).

Due to lack of work on critical sets with fractal feature, it is natural to ask how to characterize Whitney sets geometrically. Whitney posted this problem in his original paper \cite{Whi2}. The problem can be stated as follows.

{\it Given a function $f$, how far from rectifiable must a closed connected set be to be a critical set for $f$ on which $f$ is not constant?} 

\bD\label{monotone}
{\rm (\cite{WX2}) Let $\Lm$ be an arc, a homeomorphic image of the interval $[0,1]$ in $\RR^N$, and let $\Lm$ be a Whitney set. Then $\Lm$ is said to be a {\it monotone Whitney arc} if there is an increasing Whitney function $f$ on $\Lm$.}
\eD 

Xi and Wu (\cite{XW}) in 2003 gave an interesting example of a Whitney arc which is not a monotone Whitney arc.

Xi and Wu's Whitney arc $\gm$ mentioned above is not a strict Whitney set because it contains small line segments, and each Whitney function on $\gm$ must be constant on those line segments. It is not known whether there exists a strict Whitney arc which is not a monotone Whitney arc.

A mapping $f: (A, d_A)\to (B, d_B)$ between two metric spaces is said to be {\it non-expanding} if $d_B(f(x), f(y)) \le d_A(x, y)$ for $x,y \in A$.

Wen and Xi obtained the following geometric characterization of Whitney sets (see \cite[Theorem 1]{WX2}).

{\it A compact connected metric space $A$  is a Whitney set \iiff there is a non-expanding mapping from $A$ onto a monotone Whitney arc.} 

The ``if'' part is immediate because the ``pull-back'' of a Whitney function on the monotone Whitney arc is a Whitney function on $A$. The ``only if'' part can be seen as follows. 

Suppose that $f$ is a Whitney function on $A$ with $f(A)=[0,1]$. Let $B=[0,1]$. For $0\le s<t\le 1$, set
\ba{d_B(s,t)&=\inf\{\sum_{i=0}^{n-1}D(t_i,t_{i+1}): s=t_0<t_1<\dots< t_n=t\},\\
D(s,t)&=\inf _{f(x)=s,f(y)=t}d_A(x,y).}
Then $d_B$ extends in an obvious way to be a distance function on $B$. The metric space $(B, d_B)$ is a monotone Whitney arc since the identity map $\tau: B\to [0,1]$ is a monotone Whitney function on $B$. Moreover, $f: A\to B$ is a non-expanding map. For details, see \cite[p.~315]{WX2}.

\bD \label{doqa}
{\rm Let $\Lm$ be an arc, a homeomorphic image of the interval $[0,1]$ in $\RR^N$, and let $t\ge 1$. The arc $\Lm$ is said to be a {\it $t$-quasi-arc}, if there is a constant $\lm>0$ such that
\baa{\label{2017} |\Lm(x,y)|^t \le \lm|x-y|}
for each pair of points $x,y \in \Lm$, where $|\Lm(x,y)|$ is the diameter of the subarc $\Lm(x,y)$ lying between $x$ and $y$. A $1$-quasi-arc is called a {\it quasi-arc}.}
\eD

Note that (\ref{2017}) does not hold when $t<1$, because $|\Lm(x,y)|\ge|x-y|$. One can see that if an arc $\Lm$ is a $t_0$-quasi-arc, then $\Lm$ will be a $t$-quasi-arc for all $t\ge t_0$. Therefore,  each quasi-arc is a $t$-quasi-arc for each $t\ge 1$.

With the above definition of $t$-quasi-arcs, Norton (see \cite{No}) obtained the following sufficient condition for an arc $\Lm$ to be a Whitney set: if $\Lm$ is a $t$-quasi-arc and if $t$ is less than the Hausdorff dimension $\dim_H(\Lm)$ of $\Lm$, then $\Lm$ is a Whitney set.

Seeking for  necessary conditions for a $t$-quasi-arc to be a Whitney set, Norton  posed the following question (see \cite{No}): is there an arc $\Lm$ and a $C^1$ function $f$ critical but not constant on $\Lm$ such that for every subarc $\eta$ of $\Lm$ on which $f$ is not constant, $\eta$ is a $t$-quasi-arc for no $t\in [1,\infty)$?

In \cite{WX}, Wen and Xi gave an affirmative answer to the above question. They gave a Whitney function $f$ on a self-similar arc $\Lm$ such that each subarc of $\Lm$ is a $t$-quasi-arc for no $t\in [1,\infty)$. 

In Wen and Xi's work, the function $f$ is constant on some subarcs of $\Lm$, which means that $f$ is not a strict Whitney function on $\Lm$. We are interested in finding a strict Whitney function on $\Lm$. 

In \cite{No}, Norton also considered the criticality of Whitney sets. 
%\ba{\Cr(\gm)=\sup\{\beta :\exists\; g \in C^{\beta}\text{ critical but not  constant along }\gm \}.}
%For the criticality of general sets in metric space, we have the following description. 

\bD\label{cri} {\rm For a Whitney set $E$, the {\it Criticality}  of $E$ is defined to be
\ba{
\Cr(E) = \sup\{r : &\text{ there exists a non-constant function $f: E \to \RR$} \;\;\text{and
an $M>0$}\\ &\text{ such that \ $|f(x)-f(y)| \le M |x-y|^r$ \ $\forall\, x,y \in E$}\}.} }
\eD 

If $E$ is a Whitney set, then $1\le \Cr(E)\le \dim_H(E)$ (see \cite{No}). Recall that $\dim_H(E)$ is the Hausdorff dimension of $E$.

Wen and Xi worked on self-similar arcs in \cite{WX}, and obtained that each self-similar arc of Hausdorff dimension greater than $1$ is a Whitney set. 
In this paper, we obtain the following result. 
\bT\label{thm1}
Let $\Lm$ be a self-similar arc of Hausdorff dimension $s>1$. Then $\Lm$ is a strict Whitney set with criticality $\Cr(\Lm)=s$.
\eT

Theorem~\ref{thm1} improves the main result in \cite{WX} in two aspects. First, the constructed Whitney function is strictly monotone. Second, the involved H\"older component $\tilde s$ is arbitrarily close to the Hausdorff dimension $s$, hence it determines the criticality to be exactly $s$.

In Section 4, we define ``Condition $W_p$'' for a self-similar arc at the $p$-th vertex, and prove that for a  self-similar arc $\Lm$ with $\ell+1$ vertices, the Hausdorff measure function is a Whitney function on $\Lm$ \iiff Condition $W_p$ is satisfied for $p=1,\dots, \ell-1$.  We also define ``Condition $Q_p^t$'', and prove that a  self-similar arc is a $t$-quasi-arc \iiff Condition $Q_p^t$ is satisfied for all inner vertices.

In order to have a better understanding of self-similar arcs, we introduce the notion of regular self-similar arcs in Section 5. Roughly speaking, a regular self-similar arc is a self-similar arc in $\RR^2$ generated by a ``basic figure" with certain properties. One classical example of regular self-similar arc is the Koch curve.

In Section 6, we further analyze Conditions $W_p$ and $Q^t_p$ for  regular self-similar arcs, and reduce them to certain inequalities. We first prove that if the $p$-th corner angle $\theta_p>0$ then Conditions $W_p$ and $Q^t_p$ (for each $t\ge 1$) are satisfied. Consequently, a regular self-similar arc with positive corner angles is necessarily a quasi-arc and its Hausdorff measure function is a strictly monotone Whitney function. 

In case a corner angle is zero, Condition $W_p$ and Condition $Q_p^t$ are reduced to inequalities about specific parameters of the self-similar arc. By using these algebraic expressions, we could easily recognize $t$-quasi-arcs among regular self-similar arcs and determine whether the Hausdorff measure function on a regular self-similar arc is a Whitney function. 

In the last section, we provide an example of a one-parameter family of regular self-similar arcs with various features. For some value of the parameter $\tau$, the Hausdorff measure function on the self-similar arc is a strict Whitney function on the arc, and hence the self-similar arc is an $s$-quasi-arc, where $s$ is the Hausdorff dimension of the arc.  For each $t_0\ge 1$, there is a value of $\tau$ such that the corresponding self-similar arc is a $t$-quasi-arc for each $t>t_0$, but it is not a $t_0$-quasi-arc. For each $t_0>1$, there is a value of $\tau$ such that the corresponding self-similar arc is a $t_0$-quasi-arc, but it is a $t$-quasi-arc for no $t\in [1, t_0)$.

In the construction of the above mentioned one-parameter family of self-similar arcs, a crucial step in the reasoning is that the self-similar arc is a $t$-quasi-arc \iiff the parameter $\tau$ has approximation property $J_{(t-1)\log(15/7)}$. See Definition~\ref{71} for the definition of approximation property $J_a$.

The significance of the given family of self-similar arcs lies in that it provides a method to produce various examples.

\section{Self-similar arcs}
A mapping $F:\RR^N\to \RR^N$ is said to be a contractive mapping if there exists $k\in [0,1)$ such that $|F(x)-F(y)|\le k|x-y|$ for all $x,y\in \RR^N$.

A compact set $K\sbs \RR^N$ is said to be {\it invariant} with respect to a finite set $\rk S=\{S_1,\dots,S_\ell\}$ of contractive mappings on $K$, if  
\ba{K=\cup^\ell_{j=1}S_j(K).}
 In \cite{Hu}, Hutchinson gave the following theorem.
\bT Let $X=(X,d)$ be a complete metric space and let $\mathcal S=\{S_1,\dots,S_\ell\}$ be a finite set of contractive mappings on $X$. Then there exists a unique closed bounded set $K$ such that $K=\cup^\ell_{j=1}S_j(K)$. Furthermore, $K$ is compact and is the closure of the set of fixed points $s_{j_1\dots j_p}$ of finite compositions $S_{j_1}\circ\cdots\circ S_{j_\ell}$ of members of $\mathcal S$.
\eT 

A mapping $S:\RR^N\mapsto \RR^N$ is called a similitude if there is an  $r>0$ such that
\ba{|S(x)-S(y)|=r|x-y|,\;\;\text{for}\,\, x,\;y \in \RR^N.}
If $0<r<1$, we say that $S$ is a contractive similitude.

Suppose that $\rk S:=\{S_1,\dots ,S_\ell\}$ is a family of contractive similitudes with ratios $r_1,\dots,r_\ell$. Then there is a unique set $E$ satisfying
\ba{E=\cup_{i=1}^{\ell}S_{i}(E).}
The set $E$ is called the  self-similar set associated to $\rk{S}$.  

\bD\label{def2.2}
The compact self-similar set $\Lm$ associated to a family of contractive similitudes $\rk{S}=\{S_i \}_{1\leq i \leq \ell}$ is called a {\it self-similar arc} if the following two conditions are satisfied:

(1) $S_i(\Lm)\cap S_j(\Lm)$ is a singleton for $|i-j|=1$; 

(2) $S_i(\Lm)\cap S_j(\Lm)=\emptyset$ for $|i-j|>1$.
\eD
Let $\Lm$ be the self-similar arc associate to a family $\rk S=\{S_1,\dots, S_\ell\}$. The Hausdorff dimension $s$ of $\Lm$ is determined by the equation
\ba{\sum^\ell_{j=1}r_j^s=1,}
where $\{r_j\}^\ell_{j=1}$ are the contractive ratios of $\{S_j\}^\ell_{j=1}$ (see \cite{Hu}). 
We say that a self-similar arc $\Lm$ is {\it non-trivial} if the Hausdorff dimension of $\Lm$ is $s>1$, {\it i.e.}, $\Lm$ is not a line segment.

Suppose that the non-trivial self-similar arc $\Lm$ is defined by a homeomorphism $h:[0,1]\to \Lm$ so that $h(0)\in S_1(\Lm)$. For $x,y\in \Lm$, we say that $x$ precedes $y$, and write $x\prec y$, if $h^{-1}(x)<h^{-1}(y)$. Then we define intervals on $\Lm$, $[x,y]=\{z\in \Lm: x\preceq z\preceq y\}$. Now on $\Lm$, there are points $z_0\prec z_1\prec\cdots\prec z_\ell$ so that $S_j(\Lm)=[z_{j-1},z_j]$. Set $\rk S^{(k)}=\{S_{j_1}\cdots S_{j_k}:1\leq j_1,\dots,j_k\leq \ell\}$ for $k\ge 1$ and $\rk S^{(0)}=\{\Id\}$. Here $S_jS_i=S_j\circ S_i$, $S_j^k=S_j\cdots S_j$ ($k$ times), etc. 

By Definition~\ref{doqa}, for $x, y\in \Lm$ , $\Lm(x,y)$ is the subarc between $x$  and $y$. Here, we denote by $[x,y]$ the subarc from $x$ to $y$. So $\Lm(x,y)=[x,y]$ or $\Lm(x,y)=[y,x]$.

The sets $S_{j_1}\cdots S_{j_k}(\Lm)$ are intervals on $\Lm$ overlapping only at end points. Thus there are points $z_j^{(k)},j=1,\dots, \ell^k$, and a numbering $\{S_j^{(k)}:1\leq j\leq \ell^k\}$ of elements of $\rk S^{(k)}$ such that $z_0^{(k)}\prec z_1^{(k)}\prec\cdots\prec z_{\ell^k}^{(k)}$ and $S_j^{(k)}(\Lm)=[z_{j-1}^{(k)},z_j^{(k)}]$. In other works, $S_j^{(k)}$ is the unique member of $\rk S^{(k)}$ which maps $\Lm$ to $[z_{j-1}^{(k)},z_j^{(k)}]$. If $k\leq m$, we have 
\baa{\label{rela}z_j^{(k)}=z_{j\ell^{m-k}}^{(m)}.}
Note that $S_1^{(k)}$ is not necessarily equal to $S_1^k$,  because $S_1$ may be ``order reversing''.

A similitude $S_{j_1}\cdots S_{j_k}$ is order-preserving if $S_{j_1}\cdots S_{j_k}(z_0)\prec S_{j_1}\cdots S_{j_k}(z_\ell)$; otherwise it is called order-reversing.

Let $\tau$ be the function on the collection of finite sequences $(j_1,\dots,j_k)$ of members of $\{1,\dots,\ell\}$ defined by 
\begin{eqnarray}\tau(j_1,\dots,j_k)=
\begin{cases}
1, &\text{if}\;S_{j_1}\cdots S_{j_k}\;\;\text{ order-preserving,}\cr 
-1, &\text{if}\;S_{j_1}\cdots S_{j_k}\;\;\text{ order-reversing.} \cr 

\end{cases}
\end{eqnarray}
The mapping $S_{j_1}S_{j_2}$ is order-preserving if $S_{j_1}$ and $S_{j_2}$ are both order-preserving or both order-reversing; it follows that $\tau(j_1,\dots,j_k)=\tau(j_1)\cdots\tau(j_k)$.

It would be more convenient for us if $ S_1$ and $S_\ell$ are order-preserving. Of course, that is not the case in general. One might hope that when $ S_1$ and $S_\ell$ are not both order-preserving, $S_1^{(k)}$ and $S_{\ell^k}^{(k)}$ could be made order-preserving by choosing $k$ suitably. Unfortunately, that could not be achieved either. In other words, \cite[Lemma~1]{WX} is incorrect.

\bE {\rm Suppose that $S_1$ is order-preserving, $S_\ell$ is order-reversing, and $k>1$. Since $S_{\ell^m}^{(k)}$ is the unique composition $S_{j_1}\cdots S_{j_k}$ such that $S_{j_1}\cdots S_{j_k}(\Lm)\ni z_\ell$, and since $S_\ell(S_1^{k-1})(z_0)=z_\ell$, it follows that $S_{\ell^k}^{(k)}=S_\ell S_1^{k-1}$. Therefore $S_{\ell^k}^{(k)}$ is order-reversing for each $k>1$.\epf}\eE
By arguments similar to above, we obtain
\baaa{&S_1^{(k)}=S_1^k,\;\;S_{\ell^k}^{(k)}=S_\ell^k,\;\;\text{if}\;\; \tau(1)=1,\;\tau(\ell)=1,\\
&S_1^{(k)}=S_1^k,\;\;S_{\ell^k}^{(k)}=S_\ell S_1^{k-1},\;\;\text{if}\;\; \tau(1)=1,\;\tau(\ell)=-1,\\
&S_1^{(k)}=S_1S_\ell^{k-1},\;\;S_{\ell^k}^{(k)}=S_\ell^{k},\;\;\text{if}\;\; \tau(1)=-1,\;\tau(\ell)=1,\\
&S_1^{(2k)}=(S_1S_\ell)^{k},\;\;S_{\ell^k}^{(2k)}=(S_\ell S_1)^{k},\;\;\text{if}\;\; \tau(1)=-1,\;\tau(\ell)=-1,\\
&S_1^{(2k+1)}=(S_1S_\ell)^{k}S_1,\;\;S_{\ell^k}^{(2k+1)}=(S_\ell S_1)^{k}S_\ell,\;\;\text{if}\;\; \tau(1)=-1,\;\tau(\ell)=-1.}

We now define a homeomorphism $g$ from $[0,1]$ onto $\Lm$, which has properties necessary for the proofs of several theorems. Set
\ba{\Gm&=\{z_j^{(k)}:k\geq 1,\;0\leq j\leq \ell^k \},\\ Q&=\{\frac{j}{\ell^k}:k\geq 1,\;0\leq j\leq \ell^k \}\sbs [0,1].}
Let $g:Q\to \Gm$ be defined by 
\baa{\label{abcd} g(\frac{j}{\ell^k})=z_j^{(k)}.}
By (\ref{rela}), $g$ is well defined. By its very definition, the function $g$ is bijective and order-preserving, {\it i.e.,} $u<v$ implies that $g(u)\prec g(v)$. Since $Q$ is dense in $[0,1]$ and $\Gm$ is dense in $\Lm$, $g$ extends to be a homeomorphism from  $[0,1]$ onto $\Lm$.
 
Suppose that $\al\in [0,1]\setminus Q$. Then we can uniquely split $[0,1]$ into $(A,B)$ such that $a<\al<b$ for $a\in A$ and $b\in B$. Then the unique point on $\Lm$ which split $\Lm$ in $(g(A),g(B))$ is denoted by $g(\al)$. We have proved the following lemma.

\bL\label{lm2.41} There is an order-preserving homeomorphism $g:[0,1]\to \Lm$ such that for each $k\geq 1$ and each $j=1,2,\dots,\ell^k$, we have
\ba{g([\frac{j-1}{\ell^k},\frac{j}{\ell^k}])=[z_{j-1}^{(k)}, z_j^{(k)}].\;\;\epf}
\eL

%For each pair of non-negative integers $k,j$, there is a member $S_{j_1}\cdots S_{j_k}$ of $\rk S^{(k)}$ such that $g([\frac{j-1}{n^k},\frac{j}{n^k}])=S_{j_1}\cdots S_{j_k}(\Lm)$. The following lemma gives the relation between $k,j$ and $j_1,\dots,j_k$.

%\bL\label{gtlm}
%Suppose that $u\in Q$ has expansion $u=\sum_{k=1}^m u_k /n^{k}$ where $u_k=0,1,\dots, n-1$. %Then
%\ba{g(u)=S_{j_1}\cdots S_{j_m}(z_{\iota(j_1,\dots,j_m)}).}
%Where 
%$j_1=u_1+1$, $j_k=\prod_{l=1}^{k-1} \tau(j_l)(u_k-\frac{n-1}{2})+\frac{n+1}{2}$,
%$\iota(j)=(n/2)(1-\tau(j))$, and $\iota(j_1,\dots j_m)=(n/2)(1-\tau(j_1)\cdots\tau(j_m))$.
%\eL

\section{Proof of theorem \ref{thm1}}

\bL\label{measure} Let $\Lm$ be the self-similar arc associated to similitudes $S_1,\dots, S_\ell$, let $s>1$ be the Hausdorff dimension of $\Lm$, let $\tilde{s}\in (1,s)$, and let $\{\eps_k\}$ be a sequence of positive numbers with $\eps_1\le 1$ and $\eps_k\searrow 0$. Suppose that the ratios $r_j$ of $S_j$ satisfy
\baa{\label{forget}r_1^{\tilde{s}}+r_\ell^{\tilde{s}}< \sum_{j=1}^\ell r_j^{\tilde{s}}-1.}
Then there exists a number $s'\in (\tilde{s},s)$, a sequence $\{\tau_k\}$ of positive numbers with $\tau_k\searrow 0$, and a probability measure $\mu$ on $\Lm$ such that

(i) $\tau_k\leq \min(r_1^s, r_\ell^s,\eps_k^{\tilde s})$;

(ii) $\mu([x,y])>0$ \text{if} $x,y\in \Lm$ and $x\prec y$;

(iii) $\mu(S_{j_1}\cdots S_{j_k}(E))\leq (1+\eps_1)(r_{j_1}\cdots r_{j_k})^{s'}\mu(E)$ for each Borel subset $E$ of $\Lm$;

(iv) $\mu(S_{i_1}\cdots S_{i_k}(\Lm))=\tau_1\cdots \tau_k$ \text{if} $i_n=1$ \text{or} $\ell$ \text{for} $n=1,\dots,k$. 
\eL

\bpp Since $\sum_{j=2}^{\ell-1} r_j^{\tilde{s}}>1$ and $\sum_{j=2}^{\ell-1} r_j^s<1$, we see that there is an $s'\in (\tilde{s},s)$ such that $\sum_{j=2}^{\ell-1} r_j^{s'}=1$.

Let $r=\min_{2\leq j \leq \ell-1} r_j$. Choose $\gm>0$ so that $s'+\gm <s$ and $r^{-(\gm \pi^2)/6}<1+\eps_1$. Let $\tau_k'$ be such that $2\tau_k'+\sum_{j=2}^{\ell-1} r_j^{s'+\gm/k^2}=1$. Then $0<2\tau_k'<r_1^s+r_\ell^s$, and $\tau'_k\searrow 0$ as $k\to \infty$. Let $\tau_k=\min(r_1^s,r_\ell^s,\eps_k^{\tilde s},\tau_k')$. Choose $s_k$ so that $2\tau_k+\sum_{j=2}^{\ell-1} r_j^{s_k}=1$. Then $s'<s_k<s'+\gm/k^2$. Since $\tau_k\searrow 0$, we see that $s_k\searrow s'$.

For $j=1,\dots,\ell$, $k=1,2,\dots$, we define numbers $r_{jk}$ by 
\ba{r_{jk}=
\begin{cases}
r_j^{s_k}, &\text{if}\;j \neq 1,\ell,\cr 
\tau_k, &\text{if}\;j= 1\; \text{or} \;\ell. \cr 
\end{cases}}
Then we have
\baa{ \label{new1}\sum_{j=1}^\ell r_{jk}=1.}
We now define a probability measure $\mu$ by
\baa{\label{new2} \mu(\Lm)=1,\;\;\mu(S_{j_1}\cdots S_{j_k}(\Lm))=r_{j_1,1}\cdots r_{j_k,k}.} 
Equality (\ref{new1}) implies that for $k=1,2,\dots$,
\ba{\mu(S_{j_1}\cdots S_{j_{k-1}}(\Lm))=\sum_{j=1}^\ell \mu(S_{j_1}\cdots S_{j_{k-1}}S_j(\Lm)).}
Thus the definition (\ref{new2}) is consistent.

Now (i), (ii), and (iv) are satisfied; it remains to prove (iii). It suffices to show that (iii) holds for $E=S_{i_1}\cdots S_{i_n}(\Lm)$, {\it i.e.}, for arbitrary $j_1,\dots,j_k$, $i_1,\dots,i_n$, 
\ba{\mu(S_{j_1}\cdots S_{j_k}S_{i_1}\cdots S_{i_n}(\Lm))\leq (1+\eps_1)(r_{j_1}\cdots r_{j_k})^{s'}\mu(S_{i_1}\cdots S_{i_n}(\Lm)).}
We know that if $2\leq j\leq \ell-1$, then 
\ba{\frac{r_{j,k+m}}{r_{jm}}=\frac{r^{s_{k+m}}_j}{r^{s_m}_j}\leq r^{-(s_m-s_{k+m})}\leq r^{-(s_m-s')}<r^{-\gm/m^2};}
if $j=1$ or $\ell$, then 
\ba{\frac{r_{j, k+m}}{r_{jm}}={\tau_{k+m}\over\tau_k}\leq 1<r^{-\gm/m^2}.}
Thus 
\ba{\frac{r_{j, k+m}}{r_{jm}}<r^{-\gm/m^2},\;\; j=1,2,\dots,\ell,\;\; m=1,2,\cdots.}
On the other hand we have that $r_{jk}< r^{s'}_j$ for $j=1,\dots,\ell$, $k=1,2,\cdots$.
Therefore,  
\ba{\frac{\mu(S_{j_1}\cdots S_{j_k}S_{i_1}\cdots S_{i_{n}}(\Lm))}{\mu(S_{i_1}\cdots S_{i_{n}}(\Lm))}&=r_{j_1,1}\cdots r_{j_k,k} \frac{r_{i_1,(k+1)}}{r_{i_1,1}}\cdots \frac{r_{i_{\ell},(k+n)}}{r_{i_{n},n}}\\
&\leq r_{j_1}^{s'}\cdots r_{j_k}^{s'}\prod_{m=1}^{\infty}r^{-\gm/m^2}\\
&=(r_{j_1}\cdots r_{j_k})^{s'}r^{-\gm\pi^2/6}\\
&<(1+\eps_1)(r_{j_1}\cdots r_{j_k})^{s'}.}
\epp

%Now we prove Theorem \ref{thm1}.

\noindent {\it Proof of  Theorem~\ref{thm1}.\ }  Let $\tilde s\in (1, s)$ be given. We prove that there exists a function $f$ on $\Lm$, constant on no non-empty relatively open subsets of $\Lm$, and a constant $C$ such that $|f(x)-f(y)|\le C|x-y|^{\tilde s}$ for all $x, y\in \Lm$.

Suppose that $\Lm$ is the self-similar arc associated to  a family $\rk S:=\{S_1,\dots ,S_\ell\}$ of contractive similitudes with ratios $r_1,\dots,r_\ell$.   Let $g: [0,1]\to \Lm$ be the homeomorphism defined in Lemma~\ref{lm2.41}. 

Suppose that $S_j(\Lm)=[z_{j-1},z_j]$, $j=1,\cdots, \ell$. 
For each $j=1,\dots,\ell-1$, we consider sequences of points $\{\al_{jk}\}_{k=1}^{\infty}$ and $\{\beta_{jk}\}_{k=1}^{\infty}$ in $\Lm$ which are converging to $z_j$, where  
\ba{\al_{jk}=g(\frac j\ell-\frac1{\ell^k}),\;\;\;\beta_{jk}=g(\frac j\ell+\frac1{\ell^k}).}
So $\al_{j1}=z_{j-1}$\,\text{and}\ $\beta_{j1}=z_{j+1}$. 

In the following, let $\dis(X,Y)$ denote the euclidean distance between the two sets $X$, $Y$.
Set
\ba{\eps_{kj}=\min\{1,\dis([z_{j-1}, \al_{j,k+1}], [z_j, z_{j+1}]), \dis([z_{j-1}, z_j], [\beta_{j,k+1}, z_{j+1}])\},}
and $\eps_k=\min\{\eps_{kj}:j=1,\dots, \ell-1\}$. 

Since $\tilde s<s$, we have $\sum_{j=1}^\ell r_j^{\tilde s}-1>0$, hence (\ref{forget}) holds provided that $r_1$, $r_\ell$ are sufficiently small. Note that the quantity on the right side of (\ref{forget}) becomes larger when ${\mathcal S}$ is replaced by ${\mathcal S}^{(k)}$. Therefore, replacing $\mathcal S$ by ${\mathcal S}^{(k)}$ if necessary, we assume that $r_1$ and $r_\ell$ are so small that (\ref{forget}) holds.

By Lemma~\ref{measure}, there is a probability measure $\mu$ on $\Lm$ with properties (i)-(iv) specified in the lemma. 

Now we define a function $f:\Lm\to \RR$ by letting $f(x)=\mu([z_0,x])$.
By (ii), if $x\prec y$, then $f(y)-f(x)=\mu[x, y]$, hence $f$ is non-constant on each subarc of $\Lm$. We shall show that there is a constant $C>0$ such that
\baa{\label{goal} |f(x)-f(y)|\le C|x-y|^{\tilde s}, \;\;\;\text{for \ }x,y\in\Lm.}

Consider two distinct points  $x$, $y$  in $\Lm$. Let $L$ be the diameter of $\Lm$ and let $R=\max r_j$. Set 
\ba{W(x,y)=\{\kappa: \kappa\ge 0, x,y\in S^{(\kappa)}_j(\Lm)\;\;\text{for some } j,\;1\le j\le \ell^{\kappa}\}.}
Then $0\in W(x,y)$ and $W(x,y)\not=\emptyset$. When $\kappa >\log(|x-y|/L)/\log R$, the diameter of $S_j^{(\kappa)}(\Lm)$ has estimate $\diam(S_j^{(\kappa)}(\Lm))\le R^{\kappa}L<|x-y|$, which implies that $\{x,y\}\not\subset S_j^{(\kappa)}(\Lm)$, and hence $\kappa \notin W(x,y)$. Thus $\kappa \le \log(|x-y|/L)/\log R$ for $\kappa \in W(x,y)$. Let $k=\max W(x,y)$. Then $x, y\in S^{(k)}_j(\Lm)$ for some $j$ with $1\le j\le \ell^k$. Let $x', y'\in \Lm$ be such that $x=S^{(k)}_j(x')$, $y=S^{(k)}_j(y')$. Without loss of generality, we assume that $x'\prec y'$. 
Choose integers $d_1$, $d_2$ so that $x'\in S_{d_1}(\Lm)$ and $y'\in S_{d_2}(\Lm)$.
By the maximality of $k$, $d_1<d_2$.
We consider the following two cases.

Case 1. $S_{d_1}(\Lm)\cap S_{d_2}(\Lm)=\emptyset$.

By the definition of $f$ and Lemma \ref{measure}, there exists $s'\in (\tilde{s},s)$ such that 
\ba{|f(x)-f(y)|\leq \mu(S_{j_1}\cdots S_{j_k}(\Lm))\leq (1+\eps_1)(r_{j_1}\cdots r_{j_k})^{s'}\mu(\Lm).}
Let $\dl$ be the least distance between two disjoint subarcs $S_i(\Lm)$ and $S_j(\Lm)$ with $1\le i+1<j\le \ell$. Then
\ba{|x-y|&\geq \dis(S_{j_1}\cdots S_{j_k} S_{d_1}(\Lm)\;,\;  S_{j_1}\cdots S_{j_k} S_{d_2}(\Lm))\\
&= r_{j_1}\cdots r_{j_k}\dis (S_{d_1}(\Lm)\;,\;S_{d_2}(\Lm))\\
&\geq r_{j_1}\cdots r_{j_k}\dl.}
It follows that 
\ba{\frac{|f(x)-f(y)|}{|x-y|^{\tilde{s}}}\leq \frac{(1+\eps_1)(r_{j_1}\cdots r_{j_k})^{s'}}{(r_{j_1}\cdots r_{j_k})^{\tilde{s}}\dl^{\tilde{s}}} .}
Therefore,
\baa{\label{ieq1} \frac{|f(x)-f(y)|}{|x-y|^{\tilde{s}}}\leq \dl^{-\tilde{s}}(1+\eps_1).}

Case 2. $S_{d_1}(\Lm)\cap S_{d_2}(\Lm)\neq\emptyset$.

In this case, we assume that $d_1=p$ and $d_2=p+1$. Then $x'\in [z_{p-1}, z_p]=S_p(\Lm)$, $y'\in [z_p, z_{p+1}]=S_{p+1}(\Lm)$. By the maximality of $k$, we have $x'\prec z_p\prec y'$. Set $z=S_{j_1}\cdots S_{j_k}(z_p)$. Let $m$ be the least positive integer such that $x'\prec \al_{p,m+1}$. So $\al_{pm}\preceq x'\prec \al_{p,m+1}$. Similarly, let $q$ be the unique positive integer so that $\beta_{p,q+1}\prec y'\preceq \beta_{pq}$.
Then we have 
\ba{|x-y|&\geq \dis( S_{j_1}\cdots S_{j_k}([\al_{p1},\al_{p,m+1}]) \;,\; S_{j_1}\cdots S_{j_k}([\beta_{p,q+1},\beta_{p1}]) )\\
&= r_{j_1}\cdots r_{j_k}\dis( [\al_{p1},\al_{p,m+1}] \;,\;  [\beta_{p,q+1},\beta_{p1}])  \\
&\geq r_{j_1}\cdots r_{j_k}\max(\eps_{m}, \eps_{q}).}
We also have 
\ba{|f(x)-f(y)|&=\mu([x,z])+\mu([z,y])\\
&\leq \mu(S_{j_1}\cdots S_{j_k}([\al_{pm},z_p]))+\mu(S_{j_1}\cdots S_{j_k}([z_p,\beta_{pq}]))\\
&\leq (1+\eps_1)(r_{j_1}\cdots r_{j_k})^{s'} (\mu([\al_{pm},z_p])+\mu([z_p,\beta_{pq}]))\\
&=(1+\eps_1)(r_{j_1}\cdots r_{j_k})^{s'} (r_{i_1,1}\cdots r_{i_m,m}+r_{\iota_1,1}\cdots r_{\iota_q,q})\\
&\leq (1+\eps_1)(r_{j_1}\cdots r_{j_k})^{s'}( r_{i_m,m}+r_{\iota_q,q}).}
Here $r_{i_1,1}=r_p^{s_1}$, $r_{\iota_1,1}=r_{p+1}^{s_1}$,  $r_{i_j,j}=\tau_j$ for $j>1$, and $r_{\iota_n,n}=\tau_s$ for $n>1$. If $m>1$, then $r_{i_m,m}/\eps_m^{\tilde s}=\tau_m/\eps_m^{\tilde s}\le 1\le \eps_1^{ -s}$. Also, $r_{i_1,1}/\eps_1^{\tilde s}\le 1/\eps_1^{\tilde s}\le \eps_1^{ -s}$. In any case, $r_{i_m,m}/\eps_m^{\tilde s}\le \eps_1^{ -s}$. Similarly, $r_{\iota_q,q}/\eps_q^{\tilde s} \le \eps_1^{-s}$.
Therefore,
\baaa{\label{ieq2}\frac{|f(x)-f(y)|}{|x-y|^{\tilde{s}}}&\le \frac{(1+\eps_1)(r_{j_1}\cdots r_{j_k})^{s'}(r_{i_m,m}+r_{\iota_q,q}) }{(r_{j_1}\cdots r_{j_k})^{\tilde{s}}\max(\eps_{m}, \eps_{q})^{\tilde s}}\\
&\leq 2\eps_1^{-s}(1+\eps_1).}
It follows from (\ref{ieq1}) and (\ref{ieq2}) that 
\ba{ |f(x)-f(y)|\le C|x-y|^{\tilde s},} where $C=\max(2\eps_1^{-s},\dl^{-\tilde{s}})(1+\eps_1)$. Since $\tilde s>1$ is arbitrarily close to $s$, we see that  the self-similar arc $\Lm$ is a strict Whitney set and $\Cr(\Lm)=s$. \epf

\vspace{.2in}
\section{Localization}

In this section, we define ``Condition $W_p$'' for a self-similar arc at the $p$-th vertex, and prove that for a  self-similar arc $\Lm$ with $\ell+1$ vertices, the Hausdorff measure function is a Whitney function on $\Lm$ \iiff Condition $W_p$ is satisfied for $p=1,\dots, \ell-1$.  We also define ``Condition $Q_p^t$'', and prove that a  self-similar arc is a $t$-quasi-arc \iiff Condition $Q_p^t$ is satisfied for all inner vertices.

Suppose that $\Lm$ is the self-similar arc associated to  a family $\rk S:=\{S_1,\dots ,S_\ell\}$ of contractive similitudes with ratios $r_1,\dots,r_\ell$, and that the Hausdorff dimension of $\Lm$ is $s>1$. Recall that for  $x,y\in \Lm $ with $x\prec y$, $[x,y]$ is the subarc of $\Lm$ from $x$ to $y$. Let $H^s([x,y])$ be the $s$-dimensional Hausdorff measure of $[x,y]$, and let $f(x)=H^s([z_0,x])$, where $z_0$ is the ``initial point'' of $\Lm$.

As in the proof of Theorem~\ref{thm1}, for distinct points $x,y\in \Lm$, let $W(x,y)$ denote the set of positive integers $k$ such that $x,y\in S_j^{(k)}$ for some $j$ with $1\le j\le \ell^k$. 

\bD \label{df5.1}
{\rm Let $\Lm$ be a self-similar arc with $\ell+1$ vertices, and let $1\le p\le \ell-1$. The arc $\Lm$ is said to satisfy {\it Condition $W_p$} if 
\baa{\label{wpp} |f(x)-f(y)|=o(|x-y|) \;\;\text{ for } \; z_{p-1}\preceq x\preceq z_p \preceq y\preceq z_{p+1}}
or, equivalently, if 
for each $\eps>0$ there is a $\mu_p>0$ such that 
\ba{|f(x)-f(y)|\le \eps |x-y| \;\text{ whenever }  \;z_{p-1}\preceq x\preceq z_p \preceq y\preceq z_{p+1}\;\text{ and }\; |x-y|<\mu_p.}}
\eD

\bP \label{pro5.2}
The function $f$ has the property
\baa{\label{key} |f(x)-f(y)|=o(|x-y|),\;\;\; x,y\in\Lm}
\iiff $\Lm$ satisfies Condition $W_p$ for $p=1,\dots, \ell-1$.
\eP

\bpp The ``only if'' part is trivial. 

Suppose that $\Lm$ satisfies Condition $W_p$ for $p=1,\dots, \ell-1$. Let $\eps>0$ be given. Let $\mu_p>0$ be the associated number in Definition~\ref{df5.1}, $p=1,\dots,\ell-1$. Set $\mu=\min\{\mu_1,\dots,\mu_{\ell-1}\}$ and $r=\max\{r_1,\dots, r_\ell\}$.  Let $\dl_0$ be the least distance between two disjoint subarcs $S_i(\Lm)$ and $S_j(\Lm)$ with $2\le i+1<j\le \ell$. Let 
\ba{\dl_1&=(\frac{\eps \dl_0^s}{\beta})^{s-1},\;\; \dl_2=(\frac{\eps \mu^s}{\beta})^{s-1},\\
\dl&=\min(\dl_1,\dl_2), } where $s$ is the Hausdorff dimension of $\Lm$ and $\beta=H^s(\Lm)$.

Suppose that $x,y\in \Lm$ with $0<|x-y|<\dl$. Let $k=\max W(x,y)$. Then $x,y\in S^{(k)}_j(\Lm)$ for some $j$ with $1\le j\le \ell^k$. Let $x',y'\in \Lm$ be such that $x=S_j^{(k)}(x')$ and $y=S_j^{(k)}(y')$. Without loss of generality. We assume that $x'\prec y'$. There are integers $d_1\le d_2$ such that $x'\in S_{d_1}(\Lm)$, $y'\in S_{d_2}(\Lm)$. By the maximality of $k$, $d_1<d_2$. We consider the following three cases.

Case 1. $d_2-d_1>1$. Write $S_j^{(k)}=S_{j_1}\cdots S_{j_k}$. Then $(r_{j_1}\cdots r_{j_k})|x'-y'|=|x-y|<\dl_1$ implies that $(r_{j_1}\cdots r_{j_k})<\dl_1/\dl_0$. Thus
\ba{\frac{|f(x)-f(y)|}{|x-y|}&=(r_{j_1}\cdots r_{j_k})^{s-1} \frac{|f(x')-f(y')|}{|x'-y'|}\\
&< (\frac{\dl_1}{\dl_0})^{s-1}\frac{\beta}{\dl_0}=\eps. }

Case 2. $d_2-d_1=1$ and $|x'-y'|<\mu$. For convenience, set $p=d_1, p+1=d_2$. Since $|x'-y'|<\mu_p$, Condition~$W_p$ tells us that
\ba{\frac{|f(x)-f(y)|}{|x-y|}&=(r_{j_1}\cdots r_{j_k})^{s-1} \frac{|f(x')-f(y')|}{|x'-y'|}\\
&<(r_{j_1}\cdots r_{j_k})^{s-1}\eps\le \eps.}

Case 3. $d_2-d_1=1$ and $|x'-y'|\ge \mu$. As above, set $p=d_1, p+1=d_2$. Since $(r_{j_1}\cdots r_{j_k})=|x-y|/|x'-y'|<\dl_2/\mu$, it follows that 
\ba{\frac{|f(x)-f(y)|}{|x-y|}&=(r_{j_1}\cdots r_{j_k})^{s-1} \frac{|f(x')-f(y')|}{|x'-y'|}\\
&< (\frac{\dl_2}{\mu})^{s-1}\frac{\beta}{\mu}=\eps.}

Therefore, $ |f(x)-f(y)|<\eps|x-y|$ whenever $|x-y|<\dl$. The proof is complete.
\epp

Let $t\ge 1$. Set $L(x,y)=|\Lm(x,y)|^t/|x-y|$. Here $|\Lm(x,y)|$ is the diameter of the subarc $\Lm(x,y)$ of $\Lm$ between  $x$ and $y$. For $p=1,\dots,\ell-1$ we define Condition $Q_{p}^{t}$ as follows.

\bD\label{df6.1}
{\rm A self-similar arc $\Lm$ is said to satisfy {\it Condition $Q_{p}^{t}$}, if there is a constant $C_p>0$ such that $|\Lm(x,y)|^t\le C_p |x-y|$ when $z_{p-1}\preceq x\preceq z_p \preceq y\preceq z_{p+1}$}.
\eD

\bP\label{pro6.2} A self-similar arc $\Lm$ is a $t$-quasi-arc \iiff $\Lm$ satisfies Condition $Q_{p}^{t}$ for $p=1,\dots, \ell-1$.
\eP

\bpp 
The ``only if'' part is trivial. 

Suppose that $\Lm$ satisfies Condition $Q_{p}^{t}$ for $p=1,\dots, \ell-1$. Let $C_1,\dots, C_{\ell-1}$ be the numbers in Definition~\ref{df6.1} for the vertices $z_1,\dots, z_{\ell-1}$, and set  $C=\max\{C_1,\dots, C_{\ell-1}\}$.  Let $L=|\Lm|$ denote the diameter of $\Lm$. As in the proof of Proposition~\ref{pro5.2}, let $\dl_0$ denote the least distance between two disjoint subarcs $S_i(\Lm)$ and $S_j(\Lm)$. Set $M=\max(C, L^t/\dl_0)$.

Consider distinct points $x,y\in \Lm$. Let $k=\max W(x,y)$. Then $x,y\in S^{(k)}_j(\Lm)$ for some $j$ with $1\le j\le \ell^k$. Let $x',y'\in \Lm$ be such that $x=S_j^{(k)}(x')$ and $y=S_j^{(k)}(y')$. Without loss of generality, we assume that $x'\prec y'$. There are integers $d_1\le d_2$ such that $x'\in S_{d_1}(\Lm)$, $y'\in S_{d_2}(\Lm)$. By the maximality of $k$, $d_1<d_2$. We consider the following two cases.

Case 1. $d_2-d_1>1$. Write $S_j^{(k)}=S_{j_1}\cdots S_{j_k}$. Then $|\Lm(x,y)|^t=(r_{j_1}\cdots r_{j_k})^{t}|\Lm(x',y')|^t$ and $|x-y|=(r_{j_1}\cdots r_{j_k})|x'-y'|$. Therefore,
\ba{L(x,y)&=(r_{j_1}\cdots r_{j_k})^{t-1} L(x',y')\le (r_{j_1}\cdots r_{j_k})^{t-1}\frac{L^t}{\dl_0}\le \frac{L^t}{\dl_0}. }

Case 2. $d_2-d_1=1$. For convenience, set $p=d_1, p+1=d_2$. Since $z_{p-1}\preceq x'\preceq z_p\preceq y'\preceq z_{p+1}$, we have
\ba{L(x,y)&=(r_{j_1}\cdots r_{j_k})^{t-1} L(x',y')\le(r_{j_1}\cdots r_{j_k})^{t-1} C_p\le C.}

Therefore,  $L(x,y)\le M$ for distinct points $x,y\in \Lm$. By Definition~\ref{doqa}, $\Lm$ is a $t$-quasi-arc.
\epp

\section{Regular Self-similar arcs in $\RR^2$}

In this section, we study ``regular'' self-similar arcs. %We consider similitudes on the euclidean plane, which 
We identify the euclidean plane with the complex plane $\CC$ and consider the similitudes on $\CC$. It is an elementary fact that an  orientation preserving similitude $S$ is of the form $S(z)=az+b$, where $a,\;b\in \CC$, while an orientation reversing similitude $S$ has the form $S(z)=a\bar{z}+b$.

Let $\Om$ be a polygon formed by a sequence of successive segments in the plane. Suppose that $\Om$ has $\ell+1$ vertices $\{A_0, A_1,\dots, A_{\ell}\}$, and that the points  $A_1, A_{\ell-1}$ lie on segment $\overline{A_0A_\ell}$ and the point $A_{\ell-1}$ lies on segment $\overline{A_1A_\ell}$.  Suppose that there is a vertex $A_q$ such that all vertices of the polygon belong to the set $\Pi$, which is defined to be the union of the point $A_q$, the segment $\overline{A_0A_\ell}$, and the set $\Pi_0$, which is in turn defined to be the interior of triangle $A_0A_\ell A_q$. Let $\Pi_1$ be the closure of $\Pi_0$.
For $j=1, \dots, \ell$, there is a unique orientation preserving similitude $S_j$ such that $S_j(A_0)=A_{j-1}$ and $S_j(A_\ell)=A_j$. We assume that the similitudes $S_j$ are contractive, that the sets $S_j(\Pi_0)$ are pairwise disjoint, and that $S_j(A_q)\in \Pi_1$ for $j=1,\dots, \ell$. Finally we assume that
\baa{\label{fe1601}S_j(\Pi_1)\cap (\overline{A_qA_0}\cup\overline{A_qA_\ell})=\emptyset,\;\;\text{if }j\not=1,\ell,q,q+1.}

\begin{figure}[htbp] 
\centering\includegraphics[width=0.85\textwidth]{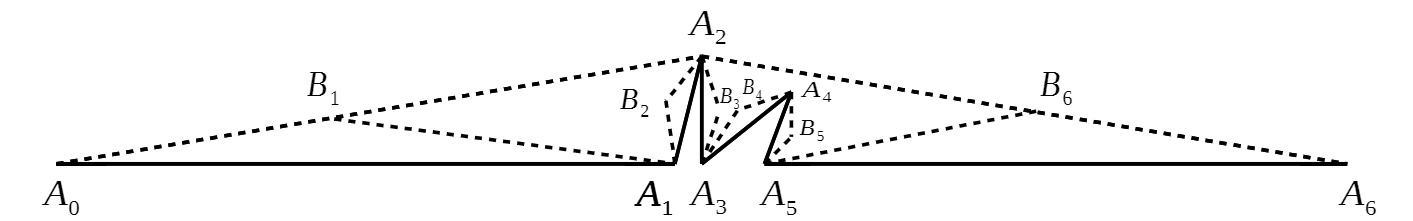} 
\caption{}\label{fig:3} 
\end{figure}

If all the above conditions are satisfied, we say that $\Om$ is a basic figure (see Figure~\ref{fig:3}), and $\Pi_1$ (and/or Triangle $A_0A_\ell A_q$, which is the union of the three sides) is the corresponding basic triangle.

Let $\Om$ be a basic figure with vertices $\{z_0,z_1,\dots ,z_{\ell}\}$ and let $\rk S=\{S_1,\dots, S_\ell\}$ be the corresponding contractive similitudes for $\Om$. Let $\Lm$ be the self-similar set associated to $\rk S$, {\it i.e.}, $\Lm$ is the unique compact set such that $\Lm=\cup_{i=1}^\ell S_i(\Lm)$. 

We now discuss under what conditions $\Lm$ is an arc. For convenience we assume that $z_0=0$, $z_\ell=1$.

\bP\label{proposition 4.1} The self-similar set $\Lm$ is an arc \iiff
\baa{\label{equ9} S_j(\Lm)\cap S_{j+1}(\Lm)=\{z_j\},\;\;j=1,\dots, \ell-1.}
\eP

\bpp Suppose that $\Lm$ is an arc. Let $1\le j\le \ell$. Then $S_j(\Lm)$ is the subarc of $\Lm$ from $z_{j-1}$ to $z_j$, and $S_{j+1}(\Lm)$ is the subarc from $z_j$ to $z_{j+1}$. Thus (\ref{equ9}) is satisfied.

Conversely, suppose that (\ref{equ9}) holds. In order to prove that $\Lm$ is an arc, we only need to prove that there is a homeomorphism between $[0,1]$ and $\Lm$.

Since $S_j$ are orientation preserving contractive similitudes for $\Om$, we know that $S_j(z)=b_jz+z_{j-1}$, where $b_j=z_j-z_{j-1}$ for $j=1,\cdots, \ell$. Now each $x\in [0,1]$ has a unique  expansion $x=\sum_{j=1}^\infty u_j /\ell^j$, where $u_j=0,\dots,\ell-1$. Recall that  
\ba{&Q:=\{\frac{j}{\ell^k}:k\geq 1,\;0\leq j\le \ell^k \}\sbs [0,1],\;\;\\
&\Gm:=\{z_j^{(k)}:k\geq 1,\;0\leq j\leq \ell^k \}\sbs\Lm.}
The function $g:Q\to \Gm$ defined in (\ref{abcd}) now has the form
\ba{g({j\over \ell^k})=z_j^{(k)}=S_j^{(k)}(1)=S_{u_1+1}S_{u_2+1}\cdots S_{u_{k-1}+1}S_{u_k+1}(1),\;\;\text{ where }\; {j\over \ell^k}=\sum_{i=1}^k {u_i\over \ell^i}.}
It is straightforward but somewhat tedious to verify that
\ba{ g(\sum_{j=1}^k \frac{u_j}{\ell^j})=\sum_{j=1}^k a_{j-1}(u_1,\cdots,u_{j-1})z_{u_j},}
where
\ba{a_0=1,\;\;a_{j}(u_1,\cdots,u_j)=\prod_{m=1}^j b_{u_m+1},\;\;j\ge 1.}

By Lemma~\ref{lm2.41}, the function $g$ extends to be  a homeomorphism $g: [0,1]\to \Lm$, which is given by
\ba{g(x)=\sum_{j=1}^\infty a_{j-1}(x)z_{u_j},}
where 
\ba{a_0(x)=1,\;\;a_{j}(x)=\prod_{m=1}^j b_{u_m+1},\;j\ge 1; \;\;x=\sum_{j=1}^\infty {u_j\over \ell^j}.}
The proposition has been proved.\epp

\begin{figure}[htbp] 
\centering\includegraphics[width=0.85\textwidth]{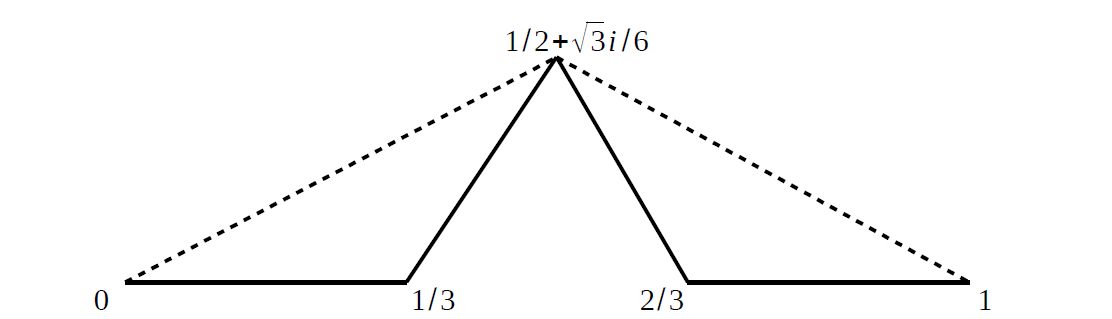} 
\caption{}\label{fig:2} 
\end{figure}

Suppose that $\Om$ is a basic figure with corresponding contractive, orientation preserving similitudes $\rk S=\{S_1,\dots, S_\ell\}$ and self-similar set $\Lm$. If $\Lm$ is an arc, then $\Lm$ is a self-similar arc by Definition~\ref{def2.2}; in this case we say that $\Lm$ is the self-similar arc generated by the basic figure $\Om$. %A {\it regular} self-similar arc is defined to be  a self-similar arc generated by some basic figure. 
For example, the Koch curve is the self-similar  arc  generated by the basic figure which is the polygon with vertices $\{0, 1/3, 1/2+\sqrt3i/6, 2/3, 1\}$ (see Figure~\ref{fig:2}), where we identify points on the complex plane with their complex number representations. Figure~\ref{fig:3} gives us an example of a basic figure with 7 vertices. Triangle $ A_2A_0A_6$ is the corresponding basic triangle.

Suppose that $\Lm$ is a self-similar arc generated by some basic figure  and the associated similitudes  $\rk S=\{S_1,\dots, S_\ell\}$ have contractive ratios $r_1,\dots, r_\ell$. The vertices of the generating basic figure are not collinear, which implies that  $r_1+\cdots +r_\ell>1$. Since the Hausdorff dimension $s$ of $\Lm$ is determined by the equation $r_1^s+\cdots +r_\ell^s=1$, it follows that $s>1$.

Let $\Om$ be a basic figure with vertices $\{z_0,z_1,\dots ,z_{\ell}\}$, and $\bigtriangleup z_q z_0 z_{\ell}$ the corresponding basic triangle. From now on we have a standing assumption that 
\ba{\Im\frac{z_q-z_0}{z_\ell-z_0}>0,}
which simplifies to $\Im z_q>0$ when $z_0=0$ and $z_\ell=1$.  For $p=1,\dots,\ell-1$, set 
\ba{\theta_p=\arg\frac{S_p(z_q)-z_p}{S_{p+1}(z_q)-z_p},\;\;0\le \theta_p<2\pi.}
Here $\theta_p$ is the argument of the fraction, so it is the angle between the two segments from the vertex $z_p$ to $S_p(z_q)$ and $S_{p+1}(z_q)$, respectively.
We call $\theta_p$ the {\it corner angle} at vertex $z_p$. Set $\theta_{\min}=\min\{\theta_1,\dots,\theta_{\ell-1}\}$. 

We now consider which points of $\Lm$ lie on the sides $\overline{z_qz_0}$ and $\ol{z_qz_\ell}$ of the basic triangle. First, points $\{S_1^j(z_q)\}$ lie on $\overline{z_qz_0}$ and accumulate at $z_0$; while $\{S_\ell^j(z_q)\}$ lie on $\overline{z_qz_\ell}$ and accumulate at $z_\ell$. For some basic figure, the side $\ol{z_qz_0}$ may contain more points of $\Lm$. For example, for the Koch curve, $q=2$, $\ell=4$, and the points $\{S_2S_4^j(z_2)\}$  lie on the side $\ol{z_2z_0}$ and accumulate at $z_2$.

For the self-similar arc $\Lm$ we define angles $\eta_1,\eta_2$ by 
\baa{\label{fe1602}\eta_1=\arg \frac{S_q(z_q)-z_q}{z_0-z_q},\;\;\eta_2=\arg\frac{z_\ell-z_q}{S_{q+1}(z_q)-z_q},\;\;0\le \eta_1,\;\eta_2<2\pi.}
(For the Koch curve $\eta_1=\eta_2=0$.) Set 
\baa{\label{fe1603}\eta_0=\min(\eta_1,\eta_2),\;\; \xi:=\theta_{\min}+\eta_0.} The angle $\xi$ is said to be the {\it characteristic angle} of $\Lm$ and of the corresponding basic figure $\Om$.

For example, in Figure~\ref{fig:3}, the basic figure $\Om$ has 7 vertices $\{A_0, A_1, A_2, A_3, A_4, A_5, A_6 \}$ with  $A_0,A_1,A_5,A_6$  collinear. We also have a family of contractive similitudes $\rk S=\{S_j: j=1,\cdots, 6\}$, where  
\ba{S_j(z)=\frac{A_j-A_{j-1}}{A_6-A_0} (z-A_0)+A_{j-1}.}
The triangle $\bigtriangleup A_2A_0A_6$ is the basic triangle, and its images under the similitudes are the smaller triangles: $\bigtriangleup B_1A_0A_1=S_1(\bigtriangleup A_2A_0A_6)$, $\bigtriangleup B_2A_1A_2=S_2(\bigtriangleup A_2A_0A_6)$, etc. Therefore, the corner angle $\theta_1=\angle B_2A_1B_1$,  $\theta_2$ is the reflex angle $\angle B_3A_2B_2$, and $\theta_4$ is the reflex angle $\angle B_5A_4B_4$, etc. The angles $\eta_1$, $\eta_2$ are $\eta_1=\angle A_0A_2B_2$, $\eta_2=\angle B_3A_2A_6$.

\bD A {\it regular} self-similar arc is a self-similar arc generated by some basic figure with a positive characteristic angle.
\eD

As in Proposition~\ref{proposition 4.1}, let $\Om$ be a basic figure with vertices $\{z_0,z_1,\dots ,z_{\ell}\}$, where $z_0=0$, $z_\ell=1$. We now express the condition (\ref{equ9}) in terms of the corner angles and other parameters of $\Om$.  Let $\Lm$ be the self-similar set generated by $\Om$. 

We fix an index $p$, where $1\le p\le \ell-1$. We first consider the case where $\theta_{p}>0$. Recall that $\Pi_1$ is the union of the basic triangle and its interior. Since $\theta_{p}>0$, we see that $S_p(\Pi_1)\cap S_{p+1}(\Pi_1)=\{z_p\}$, and $z_p=S_p(z_\ell)=S_{p+1}(z_0)$. It follows that 
\ba{S_p(\Lm)\cap S_{p+1}(\Lm)=\{z_p\}.}
Therefore condition (\ref{equ9}) holds for $i=p$ when $\theta_{p}>0$.

Now we assume that $\theta_p=0$. Let $\gamma$ be the segment $S_{p}(\ol{z_qz_\ell})$, and let $\om=S_{p+1}(\ol{z_qz_0})$. That $\theta_p=0$ means that one of the two segments $\gamma$,  $\om$ is contained in the other. Since $S_p(\Lm)\sbs S_p(\Pi_1)$ and  $S_{p+1}(\Lm)\sbs S_{p+1}(\Pi_1)$, we see  that  \ba{S_p(\Lm)\cap S_{p+1}(\Lm)\sbs S_p(\Pi_1)\cap S_{p+1}(\Pi_1)=\gm \cap \om.}  
For $j,k=0,1,2,\dots$, let
\ba{Z_j=S_pS_\ell^j(z_q),\;\;
W_k=S_{p+1}S_1^k(z_q).}
The assumption $\theta_p=0$ implies that $\eta_1,\eta_2>0$, hence $S_p(\Lm)\cap(\gm\setminus\{z_p\})=\{Z_j: j=1,2,\dots\}$.
It follows that \ba{(S_p(\Lm)\cap S_{p+1}(\Lm))\setminus \{z_p\}&=(S_p(\Lm)\cap (\gm\setminus\{z_p\}))\cap (S_{p+1}(\Lm)\cap (\om\setminus\{z_p\}))\\&=\{Z_j:j=0,1,\dots\}\cap \{W_k:k=0,1,\dots\}.} Therefore
\ba{S_p(\Lm)\cap S_{p+1}(\Lm)=\{z_p\}}\iiff
\baa{\label{szw}\{Z_j:j=0,1,2,\dots\}\cap\{W_k: k=0,1,2,\dots\}=\emptyset.}
To summarize, we conclude that 
$\Lm$ is an arc  \iiff for each $p$ with $1\le p\le \ell-1$ and $\theta_p=0$, (\ref{szw}) holds.

Since $S_j(z)=b_jz+z_{j-1}$ and $b_j=z_j-z_{j-1}$, it follows that 
\ba{Z_j=z_p-b_{p}r_{\ell}^j(1-z_q),\;\;
W_k&=z_p+b_{p+1}r_1^kz_q.}
Here $r_i=|b_i|$, for $i=1,\dots,\ell$. 
Set
\baa{ \label{feng101}\al&=|r_{p}(1-z_q)|,\;\; \beta=|r_{p+1}z_q|,\\
\label{feng102}\lm&=r_{\ell},\;\;\mu=r_1,\;\;\iota=(Z_0-z_p)/|Z_0-z_p|.}
Since $\theta_p=0$, we have
\ba{\frac{Z_j-z_p}{|Z_j-z_p|}=\frac{W_k-z_p}{|W_k-z_p|}=\iota,} 
hence
\baa{\label{feng103}&Z_j-z_p=\al\lm^j\iota,\;\;
W_k-z_p=\beta\mu^k\iota,\\
\label{feng104}&Z_j-W_k=(\al\lm^j-\beta\mu^k)\iota.}
Set 
\baa{\label{fe1304} x=-\log\lm,\;\;y=-\log\mu,\;\;u=\log(\al/\beta).} Then
\baaa{\label{feng105}u-jx+ky=\log\frac{\al\lm^j}{\beta\mu^k}.}
Therefore $Z_j\not=W_k$ for $j,k=0,1,2,\dots$ \iiff $u-jx+ky\not=0$ for $j,k=0,1,2,\dots$.

As a conclusion of the above discussion, we have the following proposition.
\bP\label{proposition4.2} Let $\Om$ be a basic figure with corner angles $\theta_p$, $p=1,\dots,\ell-1$, and let $\Lm$ be the self-similar set generated  by $\Om$. Then $\Lm$ is a regular self-similar arc \iiff for each $p$ with $1\le p\le \ell-1$ and $\theta_p=0$ the following holds:
\ba{u-jx+ky\not=0, \;\;\text{\rm for }\; j,k=0,1,2,\dots.}
\eP

\section{Reduction of Conditions $W_p$ and $Q^t_p$}

In this section we assume that $\Lm$ is the regular self-similar arc generated by  a basic figure $\Om$  and $\{S_1,\dots ,S_{\ell}\}$ are the corresponding contractive similitudes. Let $r_j$ be the ratio of $S_j$ for $j=1,\cdots, \ell$. 
Recall that  $|\Lm(x,y)|$ is  the diameter of the subarc $\Lm(x,y)$ of $\Lm$ between $x$ and $y$, and that $L=|\Lm|$. Recall also that when $x\prec y$, $[x, y]$ denotes the subarc from $x$ to $y$.

\bP \label{paiw}
Suppose that $1\le p\le \ell-1$ and that  the corner angle $\theta_p>0$.  Then $\Lm$ satisfies Condition $W_p$ .  
\eP

\bpp 
Recall that the $s$-dimensional Hausdorff measure function  $f:\Lm\to \RR$ is defined by
\ba{f(x)=H^s([z_0,x]),} 
where $z_0$ is the initial point of $\Lm$. It is clear that $f$ is non-constant on each subarc of $\Lm$. We shall  prove that there is a constant $M>0$ such that
\baa{\label{goal1} |f(x)-f(y)|\le M|x-y|^{s} \;\text{ whenever }\; z_{p-1}\prec x\prec z_p \prec y\prec z_{p+1},}
which implies (\ref{wpp}).

Suppose that  $z_{p-1}\prec x\prec z_p \prec y\prec z_{p+1}$. Let $m\ge 0$ be the greatest integer such that $x\in S_p S_\ell^m(\Lm)$. Then  $x\in S_p S_{\ell}^{m}(\Lm)$ and $x\notin S_p S_{\ell}^{m+1}(\Lm)$. Upon setting $x'=S_p S_\ell^{-m} S_p^{-1}(x)$ we obtain that $x'\in S_p(\Lm)\setminus S_pS_\ell(\Lm)$. Let $A$ denote the positive number
\ba{A:=\sup\{ |f(w)-f(z_p)|/|w-z_p|^s: w\in S_p(\Lm)\setminus S_pS_\ell(\Lm)\}.}
Since the similitude $S_p S_\ell^{-m} S_p^{-1}$ maps $x, z_p$ to $x', z_p$, respectively, it follows that
\baa{\label{feng}\frac{|f(z_p)-f(x)|}{|z_p-x|^s}=\frac{|f(z_p)-f(x')|}{|z_p-x'|^s}\le A.}
Similarly, 
\baa{\label{feng1}\frac{|f(y)-f(z_p)|}{|y-z_p|^s}\le B,}
where
\ba{B:=\sup\{ |f(w)-f(z_p)|/|w-z_p|^s: w\in S_{p+1}(\Lm)\setminus S_{p+1}S_1(\Lm)\}.}
Let 
\ba{\varrho_p:=\arg\frac{z_{p+1}-z_p}{z_{p-1}-z_p},\;\; 0<\varrho_p<2\pi,}
be the positive angle from line segment $\ol{z_{p}z_{p-1}}$ to line segment $\ol{z_{p}z_{p+1}}$. Set
\ba{\psi_p=\min(\theta_p, \varrho_p).}
It follows from the law of sines that
\baa{\label{feng2} |x-z_p|\le |x-y|\csc \psi_p,\;\;|y-z_p|\le |x-y|\csc \psi_p.}
By (\ref{feng}), (\ref{feng1}) and (\ref{feng2}), we obtain that
\ba{|f(x)-f(y)|&=|f(x)-f(z_p)|+|f(y)-f(z_p)|\\
&\le A|x-z_p|^s+B|y-z_p|^s\\
&\le (A+B)(\csc^s\psi_p)|x-y|^s.}
Thus, (\ref{goal1}) holds with $M=(A+B)(\csc^s\psi_p)$.
\epp

\bP\label{pai1q}
Suppose that $1\le p\le \ell-1$ and that  the corner angle $\theta_p>0$.  Then $\Lm$ satisfies Condition $Q_p^t$  for $t\ge 1$. 
\eP
\bpp
Fix a number $t\ge 1$. We need to prove that there is a constant $M>0$ such that 
\baa{\label{feng9}|\Lm(x,y)|^t\le M |x-y| \; \text{ whenever } \; z_{p-1}\prec x\prec z_p \prec y\prec z_{p+1}.}
Recall that $|\Lm(x,y)|$ is the diameter of the subarc $\Lm(x,y)$ of $\Lm$ between $x$ and $y$, and that $L$ is the diameter of $\Lm$.  Suppose that $x,y\in\Lm$ satisfy $z_{p-1}\prec x\prec z_p \prec y\prec z_{p+1}$. Let $m\ge 0$ be the greatest integer such that $x\in S_p S_\ell^m(\Lm)$. Then  $x\in S_p S_{\ell}^{m}(\Lm)\setminus S_p S_{\ell}^{m+1}(\Lm)$. As in the proof of Proposition~\ref{paiw} the point $x':=S_p S_\ell^{-m} S_p^{-1}(x)$ satisfies $x'\in S_p(\Lm)\setminus S_pS_\ell(\Lm)$.
Upon setting
\ba{\dl=\min(\dis(z_0, \Lm\setminus S_1(\Lm)), \dis(z_\ell, \Lm\setminus S_\ell(\Lm))),}
we obtain 
\ba{|x'-z_p|\ge \dis(z_p, S_p(\Lm)\setminus S_pS_\ell(\Lm))=r_p\dis(z_\ell, \Lm\setminus S_\ell(\Lm))\ge r_p\dl.}
Thus
\baaa{\label{feng11}\frac{|\Lm(x, z_p)|}{|x-z_p|^{1/t}}
&=r_\ell^{m(1-1/t)}\frac{|\Lm(x', z_p)|}{|x'-z_p|^{1/t}}\\
&\le \frac{|\Lm(z_{p-1}, z_p)|}{r_p^{1/t}\dl^{1/t}}\\
&= L\dl^{-1/t}r_p^{1-1/t}.}
Similarly,
\baaa{\label{feng12}\frac{|\Lm(y, z_p)|}{|y-z_p|^{1/t}}
&\le  L\dl^{-1/t}r_{p+1}^{1-1/t}.}
It follows from (\ref{feng2}), (\ref{feng11}) and (\ref{feng12}) that
\ba{|\Lm(x,y)|&\le |\Lm(x,z_p)|+|\Lm(y,z_p)|\\
&\le L\dl^{-1/t}(r_p^{1-1/t}|x-z_p|^{1/t}
+r_{p+1}^{1-1/t}|y-z_p|^{1/t})\\
&\le L(\dl^{-1}\csc \psi_p)^{1/t}(r_p^{1-1/t}
+r_{p+1}^{1-1/t})|x-y|^{1/t}.}
Therefore,
\ba{|\Lm(x,y)|^t\le  L^t\dl^{-1}(\csc \psi_p)(r_p^{1-1/t}
+r_{p+1}^{1-1/t})^t|x-y|,}
and (\ref{feng9}) has been proved.
\epp

By Propositions~\ref{pro5.2}, \ref{pro6.2}, \ref{paiw} and \ref{pai1q} we have the following
\bT\label{thm2}
Let $\Lm$ be a regular self-similar arc and let  $s=\dim_H(\Lm)$. If $\theta_{\min}:=\min\{\theta_p: p=1,\cdots, \ell-1\}>0$, then $\Lm$ is a $t$-quasi-arc for each $t\ge 1$ and the $s$-dimensional Hausdorff measure function is a Whitney function on $\Lm$. 
\eT

When the minimal corner angle $\theta_{\min}=0$, the analysis of Hausdorff measure function on $\Lm$ is more complicated. We now consider the case where $\theta_p=0$ for some $1\le p\le \ell-1$. 
As before, we assume that the three vertices of the basic triangle of  the basic figure $\Om$ under consideration  are $z_0=0$, $z_{\ell}=1$, and $z_q$ with $\Im z_q>0$.  

Let $g :[0,1]\to \Lm$ be the homeomorphism in  the proof of Proposition~\ref{proposition 4.1}. Note that $g(p/\ell)=z_p$ for $p=0,\dots,\ell$. Since $\theta_q>0$ and $\theta_p=0$, we see that $p\ne q$.

Recall that $b_p=z_{p+1}-z_p$. In section 5, we constructed two sequences of points on the self-similar arc $\Lm$, 
\ba{Z_j=S_pS_\ell^j(z_q)=z_p-b_{p}r_{\ell}^j(1-z_q),\;\;
W_k=S_{p+1}S_1^k(z_q)=z_p+b_{p+1}r_1^kz_q.} 
Set 
\baaa{\label{feng106}a=r_{p}^sH^s([z_q,z_\ell]),\;\;b=r_{p+1}^sH^s([z_0,z_q]),\;\;c=r_{p}|[z_q,z_\ell]|,\;\;d=r_{p+1}|[z_0,z_q]|,}
where $|[z_p,z_\ell]|$ denotes the diameter of the subarc of $\Lm$ from $z_p$ to $z_\ell$.
Then  we have
\ba{H^s([Z_j,z_p])=a\lm^{sj},\;\;
H^s([z_p,W_j])=b\mu^{sj},\;\; 
|[Z_j,z_p]|=c\lm^j,\;\;
|[z_p,W_j]|=d\mu^j,}
where $\lm, \mu$ are defined by (\ref{feng102}).

As usual, let $\ZZ$ denote the set of integers, let $\NN=\{1,2,\dots\}$ be the set of natural numbers, and let $\ZZ_+=\{0,1, 2,\dots,\}$.

\bL \label{feng1006} Suppose that $\Lm$ is a regular self-similar arc with $\theta_p=0$ for some $1\le p\le\ell-1$. Then there exists a constant $\Upsilon>0$ such that if $j,k\in \ZZ_+$ and if
\baaa{\label{fe1605}& x\in \Lm(Z_j,Z_{j'}),\;\;|j-j'|=1,\;\;|x-Z_j|\le |x-Z_{j'}|,\\
& y\in \Lm(W_k,W_{k'}),\;\;|k-k'|=1,\;\;|y-W_k|\le |y-W_{k'}|,}
then
 \ba{|Z_j-W_k|\le \Upsilon|x-y|.}
\eL

\bpp Since $\theta_{\min}=0$, the angles $\eta_0$, $\eta_1$ and $\eta_2$ defined by (\ref{fe1602}) and (\ref{fe1603}) are positive.  Let $\Theta$ denote the line containing the points $\{Z_j\}$ and $\{W_k\}$. Since $\eta_1, \eta_2>0$ and since (\ref{fe1601}) holds, it follows that the subarc $\Lm(Z_0, Z_1)$ intersects $\Theta$ at exactly two points $Z_0, Z_1$. The line $\Theta$ divides the plane into two half planes. Let us denote by $H_1$ the closed half plane which contains $\Lm(Z_0,Z_1)$. The other closed half plane is denoted by $H_2$. When $x\in \Lm(Z_0,Z_1)$ is sufficiently close to $Z_0$, the law of sines provides an estimate
\ba{ |z-Z_0|<(\csc \eta_0)|z-x|, \;\; z\in H_2.}
It follows that there exists a constant $C>0$ such that
\ba{|z-Z_0|\le C|z-x|,\;\;\text{ if } \; z\in H_2,\; x\in \Lm(Z_0, Z_1),\; |x-Z_0|\le |x-Z_1|.}
Similarly, there is a $C>0$ so that
\ba{|z-Z_1|\le C|z-x|,\;\;\text{ if } \; z\in H_2,\; x\in \Lm(Z_0, Z_1),\; |x-Z_0|\ge |x-Z_1|.}
It follows that for some constant $C>0$, we have
\baa{\label{fe1604}|z-Z_j| <C|z-x|}
whenever
\ba{z\in H_2,\;\;j\in \{0,1\}\;\;j'=1-j,\;\; x\in\Lm(Z_0,Z_{1}),\;\;|x-Z_j|\le |x-Z_{j'}|.}
%\le C|z-x|,\;\;\text{ if } \; z\in H_2,\; x\in \Lm(Z_0, Z_1).}
Since $\Lm(Z_j, Z_{j+1})=S_pS_\ell^jS_p^{-1}(\Lm(Z_0,Z_1))$ and since $\Theta$ and $H$ are invariant under the similitude 
$S_pS_\ell^jS_p^{-1}$, it follows that for the same constant $C$,  (\ref{fe1604}) holds 
whenever
\ba{z\in H_2,\;\;j\in \NN,\;\;|j-j'|=1,\;\; x\in\Lm(Z_j,Z_{j'}),\;\;|x-Z_j|\le |x-Z_{j'}|.}
Similarly, there exists a constant $C'>0$ such that
\ba{|z-W_k|\le C'|z-y|,}
whenever
\ba{z\in H_1,\;\;k\in \NN,\;\;|k-k'|=1,\;\; y\in\Lm(W_k,W_{k'}),\;\;|y-W_k|\le |y-W_{k'}|.}
Now suppose that $j,k\in \NN$ and (\ref{fe1605}) holds. Then
\ba{|Z_j-W_k|\le C|x-W_k|\le CC'|x-y|.}
The proof is complete.\epp

\bP\label{a0iw}
Suppose that $1\le p\le\ell-1$ and  $\theta_p=0$. Then $\Lm$ satisfies Condition $W_p$ \iiff
\baa{\label{wp} (\lm^j+\mu^k)^s=o(|\al\lm^j-\beta\mu^k|),\;\;j,k=0,1,2,\dots.}
\eP

%{The equation (\ref{wp}) means that for given $\eps>0$ there exists $\dl>0$ such that if $j,k\in\NN$ and if $|A\lm^j-B\mu^k|<\dl$ then $ (\lm^j+\mu^k)^s\le \eps |A\lm^j-B\mu^k|$.}

\bpp
Suppose that $\Lm$ satisfies Condition $W_p$. By Definition~ \ref{df5.1}, we know that (\ref{wpp}) holds for all $Z_j$ and $W_k$, {\it i.e.,}
\baa{\label{feng1001} H^s([Z_j,W_k])=|f(Z_j)-f(W_k)|=o(|Z_j-W_k|).} 
By (\ref{feng106}), we have 
\ba{H^s([Z_j,W_k])=a\lm^{sj}+b\mu^{sk},}
which, together with (\ref{feng104}) and (\ref{feng1001}), implies that
\baa{\label{feng1002}a\lm^{sj}+b\mu^{sk}=o(|\al\lm^j-\beta\mu^k|).}
The H\"older inequality tells us that  
\baa{\label{feng1004}(\lm^j+\mu^k)^s\le F(a\lm^{sj}+b\mu^{sk}),}
where $F=(a^{-1/(s-1)}+b^{-1/(s-1)})^{s-1}$. Now (\ref{wp}) is a consequence of (\ref{feng1002}) and (\ref{feng1004}).

Conversely, suppose that (\ref{wp}) holds.  
Let $x,y\in \Lm$ be such that $z_{p-1}\prec x\prec z_p\prec y\prec z_{p+1}$.
Let $m$ be the least positive integer such that $x\prec Z_m$, so $x\in[Z_{m-1},Z_m]$. If $|x-Z_{m-1}|\le |x-Z_m|$ let $j=m-1$ and $j'=m$; otherwise, let $j=m, j'=m-1$. In either case we have 
\baa{\label{fe1301}  x\in \Lm(Z_j,Z_{j'}),\;\;|x-Z_j|\le |x-Z_{j'}|,\;\; |j-j'|=1.} Similarly, there are integers $k,k'\ge 0$ such that 
\baa{\label{fe1302} y\in \Lm(W_k,W_{k'}),\;\; |y-W_k|\le |y-W_{k'}|,\;\;
|k-k'|=1.}
By Lemma~\ref{feng1006}, we have
\ba{|Z_j-W_k|\le \Upsilon|x-y|,}
which, together with (\ref{wp}), implies that
\baa{\label{feng1007} (\lm^j+\mu^k)^s=o(|x-y|).}
Setting $Z_{-1}=Z_0$ and $W_{-1}=W_0$, we have
\ba{H^s([x,y])\le &H^s([Z_{j-1},W_{k-1}])\\
\le&a\lm^{s(j-1)}+b\mu^{s(k-1)}\\
\le &(a\lm^{-s}+b\mu^{-s})(\lm^{sj}+\mu^{sk})\\
\le &(a\lm^{-s}+b\mu^{-s})(\lm^j+\mu^k)^s.}
The last inequality and (\ref{feng1007}) imply that
\ba{H^s([x,y])=o(|x-y|).}
\epp

\bP \label{a0ht}
Suppose that $t\ge 1$, $1\le p\le\ell-1$,  and  $\theta_p=0$. Then $\Lm$ satisfies Condition $Q_{p}^{t}$ \iiff there exists a constant $C>0$ such that
\baa{\label{tp} (\lm^j+\mu^k)^t \le C|\al\lm^j-\beta\mu^k|,\;\;j,k=0,1,2,\dots.}
\eP

\bpp
Suppose that $\Lm$ satisfies Condition $Q_{p}^{t}$. Then there exists a constant $C_p>0$ such that 
\baa{\label{fe1303} |[Z_j, W_k]|^t\le C_p |Z_j-W_k|,\;\; j,k=0,1,\dots.}
We have the following estimate
\ba{\min(c,d)(\lm^j+\mu^k)\le c\lm^j+d\mu^k=|[Z_j,z_p]|+|[z_p,W_k]|\le 2|[Z_j,W_k]|,}
which, together with (\ref{feng104}) and (\ref{fe1303}), implies that 
\ba{(\lm^j+\mu^k)^t \le C'|[Z_j,W_k]|^t\leq C|\al\lm^j-\beta\mu^k|,} 
where $C=C'C_p=\{2/\min(c,d)\}^t C_p$. Thus (\ref{tp}) holds.

Conversely, suppose that there exists a constant $C>0$ such that (\ref{tp}) holds. Let $x,y\in \Lm$ be such that $z_{p-1}\prec x\prec z_p\prec y\prec z_{p+1}$. We need to prove that there exists constant $M>0$ such that  
\ba{|[x,y]|^t\le M |x-y|.}
As in the proof of previous proposition, there exist integers $j,j', k,k'$ such that (\ref{fe1301}) and (\ref{fe1302}) hold.
By Lemma~\ref{feng1006}, we have
\ba{|Z_j-W_k|\le \Upsilon|x-y|,}
which, together with (\ref{feng104}) and (\ref{tp}), implies that
\baa{\label{feng1008} (\lm^j+\mu^k)^t\le C\Upsilon |x-y|.}
Now
\ba{|[x,y]|^t\le &|[Z_{j-1},W_{k-1}]|^t\\
\le & (|[Z_{j-1},z_p]|+|[z_p,W_{k-1}]|)^t\\
\le&(c\lm^{j-1}+d\mu^{k-1})^t\\
\le &(c\lm^{-1}+d\mu^{-1})^t(\lm^j+\mu^k)^t.}
The last inequality and (\ref{feng1008}) imply that
\ba{|[x,y]|^t\le M |x-y|,}
where $M=(c\lm^{-1}+d\mu^{-1})^tC\Upsilon$.
\epp

\bP\label{thm3} 
Let $\Lm$ be a regular self-similar arc and let  $s=\dim_H(\Lm)$. If the $s$-dimensional Hausdorff measure function $f$ is a Whitney function on $\Lm$, then $\Lm$ is an $s$-quasi-arc. If $\Lm$ is a $t$-quasi-arc for some $t$ with $s>t\ge 1$, then $f$ is a Whitney function on $\Lm$.
\eP

\bpp By Theorem~\ref{thm2} and Proposition~\ref{a0ht}, the self-similar arc $\Lm$ is a $t$-quasi arc \iiff for each $p$ with $\theta_p=0$, one has 
\baa{\label{wp9} (\lm^j+\mu^k)^t=O(|\al\lm^j-\beta\mu^k|),\;\text{where}\;j,k=0,1,2,\dots.} 
By Theorem~\ref{thm2} and Proposition~\ref{a0iw}, the $s$-dimensional Hausdorff measure function on $\Lm$ is a Whitney function \iiff for each $p$ with $\theta_p=0$, one has  
\baa{\label{wp1} (\lm^j+\mu^k)^s=o(|\al\lm^j-\beta\mu^k|),\;\text{where}\;j,k=0,1,2,\dots.} 
The proposition follows because (\ref{wp1}) implies (\ref{wp9}) when $t=s$, and because (\ref{wp1}) follows from (\ref{wp9}) when $1\le t<s$.\epp

The second part of Proposition~\ref{thm3} is contained in the result of  Norton  \cite{No} mentioned in the introduction of this paper.

By Proposition~\ref{a0ht}, Condition~$Q_p^t$ is reduced to an inequality (\ref{tp}). In the following proposition it is further reduced to an inequality of a certain form which is more convenient for determining whether a self-similar arc is a $t$-quasi-arc and which is directly related to the degree to which a number $u$ is approximated by numbers of the form $jx-ku$, where $x,y$ are fixed positive numbers and $j,k$ are non-negative integers.

Recall that $\lm, \mu, x, y, u$ are defined by (\ref{feng102}) and (\ref{fe1304}).

\bP\label{pro8.1} 
Suppose that $1\le p\le \ell-1$ and $\theta_p=0$. Then $\Lm$ satisfies Condition $Q_{p}^{t}$ \iiff there exists a constant $M>0$ such that
\baa{\label{tpp}e^{-j(t-1)x}\le M|u-jx+ky|,\;\;j,k=0,1,2,\dots.}
\eP

\bpp\label{pro0.18}
By Proposition~\ref{a0ht}, we only need to show that there exists a constant $C>0$ such that (\ref{tp}) holds \iiff there exists a constant $M>0$ such that (\ref{tpp}) holds.

Suppose that there exists no constant $C$ such that (\ref{tp}) holds. Then there are increasing sequences $\{j_n\}$ and $\{k_n\}$ of positive integers such that
\baa{\label{tppp}\lim_{n\to\infty} \frac{ |\al\lm^{j_n}-\beta\mu^{k_n}|}{(\lm^{j_n}+\mu^{k_n})^t}=0.}
Since $(\lm^{j_n}+\mu^{k_n}) \ge (\lm^{j_n}+\mu^{k_n})^t$ when $n$ is large enough,  we see that
\baa{\lim_{n\to\infty} \frac{ |\al\lm^{j_n}-\beta\mu^{k_n}|}{\lm^{j_n}+\mu^{k_n}} \le \lim_{n\to\infty} \frac{ |\al\lm^{j_n}-\beta\mu^{k_n}|}{(\lm^{j_n}+\mu^{k_n})^t}=0.}
This implies that when $n$ is sufficiently large, the quotient $|\al\lm^{j_n}-\beta\mu^{k_n}|/(\lm^{j_n}+\mu^{k_n})$ does not exceed $\al/2$, hence we have 
\ba{\al\lm^{j_n}&\le |\al\lm^{j_n}-\beta\mu^{k_n}|+ \beta\mu^{k_n} \le \frac{\al}{2} (\lm^{j_n}+\mu^{k_n})+ \beta\mu^{k_n},\\
\frac{\al}{2} \lm^{j_n} &\le (\frac{\al}{2}+\beta)\mu^{k_n}.}
It follows that there is a constant $Q_1$ such that $\lm^{j_n}\le Q_1^{-1}\mu^{k_n}$ for all $n$. Similarly, there is a constant $Q_2$ such that $\mu^{k_n}\le Q_2\lm^{j_n}$ for all $n$. Therefore, 
\baa{\label{tppp1} Q_1 \lm^{j_n} \le \mu^{k_n}\le Q_2 \lm^{j_n}.}
%Thus we have $\lm^{j_n}+\mu^{k_n}\le (1+Q_2)\lm^{j_n}$, which together with (\ref{tppp}) implies that \baa{\label{tppp6}\lim_{n\to\infty} \frac{ A\lm^{j_n}-B\mu^{k_n}}{\lm^{tj_n}} =0.}
Now 
\baa{\label{tppp5}\beta\frac{(\al/\beta)\lm^{j_n}\mu^{-k_n}-1}{\lm^{(t-1)j_n}}=\frac{\al\lm^{j_n}-\beta\mu^{k_n}}{(\lm^{j_n}+\mu^{k_n})^t}\Big(\frac{\lm^{j_n}+\mu^{k_n}}{\lm^{j_n}}\Big)^t\frac{\lm^{j_n}}{\mu^{k_n}}.}
The first factor on the right side of (\ref{tppp5}) tends to $0$ as $n \to \infty$
by (\ref{tppp}), while the second and third factors are bounded above and below because of (\ref{tppp1}). It follows that 
\baa{\label{tppp7}\lim_{n\to\infty} \frac{ (\al/\beta)\lm^{j_n}\mu^{-k_n}-1}{\lm^{(t-1)j_n}}=0.}
Since the denominator in (\ref{tppp7}) is $\le 1$, it follows that the numerator tends to $0$ as $n\to \infty$.
%By inequalities (\ref{tppp1}), we have
%\ba{\frac{\lm^{tj_n}}{\mu^{k_n}}=\lm^{(t-1)j_n} \frac{\lm^{j_n}}{\mu^{k_n}}\ge \frac{1}{Q_2} \lm^{(t-1)j_n},}
%hence we have the estimate 
%\baa{\lim_{n\to\infty} \frac{ (A/B)\lm^{j_n}\mu^{-k_n}-1}{\lm^{(t-1)j_n}} \le \lim_{n\to\infty} \frac{ (A/B)\lm^{j_n}\mu^{-k_n}-1}{(Q_2/B)(\lm^{tj_n}/\mu^{k_n})}=\frac{1}{Q_2} \lim_{n\to\infty} \frac{ A\lm^{j_n}-B\mu^{k_n}}{\lm^{tj_n}}=0.}
%So we obtain that
%\baa{\label{tppp4} \lim_{n\to\infty} \log((A/B)\lm^{j_n}\mu^{-k_n}) \le \lim_{n\to\infty}  (A/B)\lm^{j_n}\mu^{-k_n}-1) \le \lim_{n\to\infty} \lm^{(t-1)j_n}=0 .}
By (\ref{tppp7}) and the equality $\lim_{w\to1}[(\log w)/(w-1)]= 1 $, we have 
\baa{\label{tppp3} \lim_{n\to\infty} \frac{ \log((\al/\beta)\lm^{j_n}\mu^{-k_n})}{\lm^{(t-1)j_n}}=0.}
Substituting $x=-\log\lm$, $y=-\log\mu$, and $u=\log(\al/\beta)$ into (\ref{tppp3}) yields that
\baa{\label{tppp2} \lim_{n\to\infty}\frac{u-j_nx+k_ny}{e^{-j_n(t-1)x}}=0.}
Thus there exists no $M$ such that (\ref{tpp}) holds.  

Conversely, suppose that there exists no $M$ such that (\ref{tpp}) holds. Then there are increasing sequences $\{j_n\}$ and $\{k_n\}$ of positive integers such that (\ref{tppp2}) holds, hence the equivalent equalities (\ref{tppp3}) and (\ref{tppp7}) hold. Since the numerator in (\ref{tppp7}) tends to 0 as $n\to\infty$, it follows that $1/2<(\al/\beta)\lm^{j_n}\mu^{-k_n}<2$ for $n$ large enough, which implies (\ref{tppp1}).
Then (\ref{tppp}) follows from (\ref{tppp1}) and (\ref{tppp7}). Therefore there exists no $C$ such that (\ref{tp}) holds.
\epp

\bP\label{thmrat}
Suppose that $1\le p\le \ell-1$ and $\theta_p=0$. Then $\Lm$ satisfies Condition $Q_{p}^{1}$ \iiff $x/y$ is rational.
\eP

\bpp By Propositions~\ref{a0ht} and \ref{pro8.1}, $\Lm$ satisfies Condition $Q_{p}^{1}$ \iiff the inequality (\ref{tpp}) holds when $t=1$.

Suppose that $\tau:=x/y$ is rational. Then the set \ba{\Pi:=\{j\tau-k: j,k\in\NN\}}
is discrete. By Proposition~\ref{proposition4.2}, the distance from the point $u/y$ to $\Pi$ is positive. It follows that $|u-jx+ky|>\dl$ for some $\dl>0$. Thus the inequality  (\ref{tpp}) holds with $M=1/\dl$. 

Suppose that $\tau=x/y$ is irrational. We show that the set $\Pi$ 
is dense in $\RR$, which implies that (\ref{tpp}) does not hold with $t=1$ and that $\Lm$ does not satisfy Condition $Q_{p}^{1}$. 
Let $c\in \RR$ and $\eps>0$. There exist $j_0,k_0\in\ZZ$ such that
\baa{ \label{f101}|(j_0\tau-k_0)-c|<\eps/2.}
By Dirichlet's Approximation Theorem (see, {\it e.g.}, \cite[p.~143]{Ap}), there are positive integers $j',k'>\max(|k_0|, |j_0|, 2\tau/\eps)$ such that
\ba{|\frac{j'}{k'}-\frac1\tau|<\frac1{k'^2}<\frac\eps{2k'\tau},}
hence
\baa{\label{f102} |j'\tau-k'|<\eps/2.}
Set $j=j_0+j'$ and $k=k_0+k'$. Then $j,k\in\NN$. It follows from (\ref{f101}) and (\ref{f102}) that $|(j\tau-k)-c|<\eps$. Therefore, $\Pi$ is dense in $\RR$.\epp

\section{A one-parameter family of self-similar arcs}

In this section, we  construct and examine a one-parameter family of regular self-similar arcs with $\theta_{p}=0$ for some fixed $p$. For different values of the parameter $\tau$, the corresponding regular self-similar arcs have various features. It turns out that the self-similar arc satisfies Condition $Q_p^t$ \iiff  the number $\tau$ satisfies a certain ``approximation property'' $J_{(t-1)\zeta}$, where $\zeta=\ln(15/7)$. We now define
approximation property $J_a$, $a>0$, of irrational numbers.

\bD \label{71} {\rm Let $a>0$. An irrational number $\tau$ is said to have {\it approximation property} $J_a$ if 
\baa{\label{def}\exists\; C>0,\;\; |\tau-k/j|\ge Cj^{-1} e^{-aj},\;\;\forall\; k\in \ZZ, j\in \NN.}}
\eD

It follows directly from the definition that if $\tau$ has approximation property $J_{a_0}$ then $\tau$ has approximation property $J_a$ for each $a>a_0$. By Liouville's Approximation Theorem (see, {\it e.g.}, \cite[p.~146]{Ap}), each algebraic irrational number $\tau$ satisfies $|\tau-k/j|>Cj^{-m}$, where $m$ is the degree of the irreducible polynomial with integer coefficients of which $\tau$ is a root, hence $\tau$ has approximation property $J_a$ for each $a>0$.

\bT\label{7.11} Let $a_0>0$ and let $\nu\in \NN$. There exists a transcendental number $\tau$ with $1<\tau<1+2^{-\nu}$ such that $\tau$ has approximation property $J_{a_0}$, but $\tau$ has approximation property $J_a$ for no $a\in (0,a_0)$.
\eT

\bpp Define a number $\tau$ by 
\baaa{\label{number tau}\tau &=1+2^{-n_1}+2^{-n_2}+\cdots,\\
n_1&=m_1, n_2=m_1+m_2, n_3=m_1+m_2+m_3,\dots,\\
m_1&=\max(8,\nu+1), \;\; m_{i+1}=\max(\lceil 2^{n_i}a_0/\log 2\rceil-n_i,n_i+2) \;\text{ for }\;i\ge1,}
where $\lceil \cdot\rceil$ is the ceiling function, {\it i.e.}, $\lceil u\rceil$ is the least integer greater than or equal to  $u$. Since $m_i\to\infty$ as $i\to\infty$, we see that
\ba{\tau=1+2^{-m_1}(1+2^{-m_2}+2^{-m_2-m_3}+\cdots)<1+2^{-m_1+1}\le 1+2^{-\nu}.}

For $i\ge1$, set $j_i=2^{n_i}$ and $k_i=j_i\tau_i$, where
\baa{\label{ti}\tau_i =1+2^{-n_1}+2^{-n_2}+\cdots+2^{-n_i}.}
Then $k_i$ is an integer, and 
\ba{j_i\tau=k_i+(2^{-m_{i+1}}+2^{-m_{i+1}-m_{i+2}}+\cdots).}
It follows that
\baa{\label{ineq1}0<j_i\tau-k_i<2^{-m_{i+1}+1}.}
From the definition of $m_{i+1}$ in (\ref{number tau}), we see that there is an   $i_0\in\NN$ such that for $i\ge i_0$ we have  
\baa{\label{ineq2}2^{m_{i+1}-1}< 2^{-n_i}\exp(2^{n_i}a_0)\le 2^{m_{i+1}}.}
Combining inequalities (\ref{ineq1}) and (\ref{ineq2}), we obtain that
\ba{0<j_i\tau-k_i<2j_ie^{-a_0j_i}.}
Consider a fixed number $a\in (0,a_0)$. The above inequality tells us that for $i\ge i_0$,
\ba{\frac{|\tau-k_i/j_i|}{j_i^{-1} e^{-aj_i}}<2j_i e^{(a-a_0)j_i}.}
Since $a-a_0<0$ and hence the right side of the above inequality tends to 0 as $i$ approaches $\infty$, we see that (\ref{def}) does not hold. Thus $\tau$ does not have approximation property $J_a$.

Now we assume that $j\ge 2^{n_{i_0}}$ and $k$ is an arbitrary integer. Then there is an $i\ge i_0$ such that $2^{n_i}\le j< 2^{n_{i+1}} $. By (\ref{ineq2}), the integers $n_i$ and $n_{i+1}$ satisfy 
\baa{\label{ab1}2^{n_{i+1}-1}< \exp(2^{n_i}a_0)\le 2^{n_{i+1}}.}
In order to obtain a lower bound for $|j\tau-k|$, we write \baa{\label{twoterm}j\tau-k=(j\tau_{i+1}-k)+j(\tau-\tau_{i+1}).} Recall that $\tau_i$ is defined by (\ref{ti}). Since $2^{n_{i+1}}\tau_{i+1}$ is an odd integer, and since $j$ is not a multiple of $2^{n_{i+1}}$, we see that $2^{n_{i+1}}j\tau_{i+1}$ is not a multiple of $2^{n_{i+1}}$. It follows that
$|2^{n_{i+1}}j\tau_{i+1}- 2^{n_{i+1}}k|\ge1$, and therefore
\baa{\label{ab2}|j\tau_{i+1}- k|\ge2^{-n_{i+1}}.}
For the second term on the right side of (\ref{twoterm}) we have
\ba{j(\tau-\tau_{i+1})&<2^{n_{i+1}}(2^{-n_{i+2}}+2^{-n_{i+3}}+\cdots)\\ 
&=2^{-m_{i+2}}+2^{-m_{i+2}-m_{i+3}}+\cdots\\
&<2^{-m_{i+2}+1}.}
Since $m_{i+2}\ge n_{i+1}+2$, by the definition of $m_i$, the right side of the last inequality is $\le 2^{-n_{i+1}-1}$. Thus 
\baa{ \label{ab3}0<j(\tau-\tau_{i+1})< 2^{-n_{i+1}-1}.}
Now inequalities (\ref{twoterm}), (\ref{ab2}) and (\ref{ab3}) tell us that
\baa{\label{ab}|j\tau- k|\ge2^{-n_{i+1}-1}.}
From (\ref{ab}) and (\ref{ab1}), we obtain that
\ba{|j\tau- k|&\ge2^{-n_{i+1}-1}\\
&>(1/4) \exp(-2^{n_i}a_0)\\
&\ge (1/4) e^{-a_0j}.}
Therefore, the inequality in (\ref{def}), with $a$ replaced by $a_0$, holds with $C=1/4$ as long as $j\ge 2^{n_{i_0}}$. This implies that $\tau$ has approximation property $J_{a_0}$.

Finally, since $\tau$ does not have approximation property $J_a$ for $a<a_0$, it cannot be an algebraic number. Thus $\tau$ is a transcendental number.
\epp

\bT\label{7.12} Let $a_0>0$ and let $\nu\in \NN$. There exists a transcendental number $\tau$ with $1<\tau<1+2^{-\nu}$ such that $\tau$ has approximation property $J_{a}$ for each $a>a_0$, but $\tau$ does not have approximation property $J_{a_0}$.
\eT

\bpp Define a number $\tau$ by 
\baaa{\label{tau}\tau &=1+2^{-n_1}+2^{-n_2}+\cdots,\\
n_1&=m_1, n_2=m_1+m_2, n_3=m_1+m_2+m_3,\dots,\\
m_1&=\max(8,\nu+1), \;\; m_{i+1}=n_i+\max(\lceil 2^{n_i}a_0/\log 2\rceil,2) \;\text{ for }\;i\ge1.}
Then $\tau$ satisfies $1<\tau<1+2^{-\nu}$. As in the previous theorem, $\tau$ is a transcendental number because we shall show that $\tau$ does not have property $J_{a_0}$.

Setting  $j_i=2^{n_i}$ and $k_i=j_i\tau_i$, we obtain that
\baa{\label{np1}|j_i\tau-k_i|<2^{-m_{i+1}+1}.}
By the definition of $m_{i+1}$, there is an   $i_0\in\NN$ such that for $i\ge i_0$ we have  
\baa{\label{np2}2^{m_{i+1}-1}< 2^{n_i}\exp(2^{n_i}a_0)\le 2^{m_{i+1}}.}
We then combine (\ref{np1}) and (\ref{np2}) to obtain
\ba{|j_i\tau-k_i|<2j_i^{-1}e^{-a_0j_i},}
which implies that for $i\ge i_0$,
\ba{\frac{|\tau-k_i/j_i|}{j_i^{-1} e^{-a_0j_i}}<2j_i^{-1} .}
Thus $\tau$ does not have approximation property $J_{a_0}$.

Let $a>a_0$. We now prove that $\tau$ has approximation property $J_a$. Choose an integer $i_1\ge i_0$ such that whenever $i\ge i_1$, the following inequality holds:
\baa{ \label{np3}2^{-2 n_i-2} \exp(-2^{n_i}a_0)>  \exp(-2^{n_i}a).}
Assume that $i\ge i_1$, $2^{n_i}\le j< 2^{n_{i+1}}$, and $k$ is an arbitrary integer. Similar to the previous proof, we have
\baa{\label{np4}2^{n_{i+1}-1}< 2^{2n_i}\exp(2^{n_i}a_0)\le 2^{n_{i+1}},}
and
\baa{\label{np5}|j\tau- k|\ge2^{-n_{i+1}-1}.}
From (\ref{np3}), (\ref{np4})  and (\ref{np5}), we obtain that
\ba{|j\tau- k|&\ge2^{-n_{i+1}-1}\\
&>2^{-2 n_i-2} \exp(-2^{n_i}a_0)\\ &>  \exp(-2^{n_i}a)\\
&\ge  e^{-aj}.}
Therefore, the inequality in (\ref{def}) holds with $C=1$ as long as $j\ge 2^{n_{i_1}}$. This implies that $\tau$ has approximation property $J_{a}$.
\epp

\bT\label{7.13}
Let $\nu\in \NN$. Then there exists a transcendental number $\tau$ with $1<\tau<1+2^{-\nu}$ such that $\tau$ has approximation property $J_{a}$ for each $a>0$.
\eT

\bpp Define a number $\tau$ by 
\baaa{\label{tau3}\tau &=1+2^{-n_1}+2^{-n_2}+\cdots,\\
n_1&=m_1, n_2=m_1+m_2, n_3=m_1+m_2+m_3,\dots,\\
m_1&=\max(8,2\nu), \;\; m_{i+1}= 2^{n_i/2} \;\text{ for }\;i\ge1.}
Then $1<\tau<1+2^{-\nu}$, as in the previous theorem. It is clear that for each $i$, $m_i$ is an integer, and $m_{i+1}>n_i+2$, which will be needed later.

Setting  $j_i=2^{n_i}$ and $k_i=j_i\tau_i$, we obtain that
\baa{\label{np7}|j_i\tau-k_i|<2^{-m_{i+1}+1}=2^{-\sqrt{j_i}+1}.}
It follows that for each positive integer $n$,
\ba{\lim_{i\to\infty}\frac{|j_i\tau-k_i|}{j_i^{-n}}=0.}
By Liouville's Approximation Theorem, $\tau$ must be a transcendental number.

Let $a>0$. We now prove that $\tau$ has approximation property $J_a$. Choose $i_0$ so that when $i\ge i_0$, we have
\ba{2^{-n_{i+1}-1}=2^{-n_i-2^{n_i/2}-1}>\exp(-2^{n_i}a).}
Assume that $i\ge i_0$, $2^{n_i}\le j< 2^{n_{i+1}}$, and $k$ is an arbitrary integer. Similar to the previous proof, since $m_{i+1}>n_i+2$, we see that
\ba{|j\tau- k|\ge2^{-n_{i+1}-1}.}
It follows that
\ba{|j\tau- k|>\exp(-2^{n_i}a)\ge e^{-aj}.}
Therefore, $\tau$ has approximation property $J_a$.
\epp

\bT\label{7.14} Let $\nu\in \NN$. Then there exists a transcendental number $\tau$ with $1<\tau<1+2^{-\nu}$ such that $\tau$ has approximation property $J_{a}$ for no $a>0$.
\eT

\bpp Define a number $\tau$ by 
\baaa{\label{tau4}\tau &=1+2^{-n_1}+2^{-n_2}+\cdots,\\
n_1&=m_1, n_2=m_1+m_2, n_3=m_1+m_2+m_3,\dots,\\
m_1&=\max(8,\nu+1), \;\; m_{i+1}= 2^{2n_i} \;\text{ for }\;i\ge1.}
Then $1<\tau<1+2^{-\nu}$, as in the previous theorem. 

Setting  $j_i=2^{n_i}$ and $k_i=j_i\tau_i$, we obtain, as in the proof of the previous theorem, that
\baa{\label{np8}|j_i\tau-k_i|<2^{-m_{i+1}+1}=2^{-j_i^2+1}.}
Consider a fixed number $a>0$. Then (\ref{np8}) implies that
\ba{\lim_{i\to\infty}\frac{|j_i\tau-k_i|}{e^{-aj_i}}=0,}
and hence $\tau$ does not have approximation property $J_a$.
By Liouville's Approximation Theorem, $\tau$ is necessarily  a transcendental number.
\epp

Now we construct a one-parameter family of regular self-similar arcs. We start by constructing a family of basic figures depending on a parameter $\tau$ with $1<\tau<1.001$. 

\begin{figure}[htbp] 
\centering\includegraphics[width=0.85\textwidth]{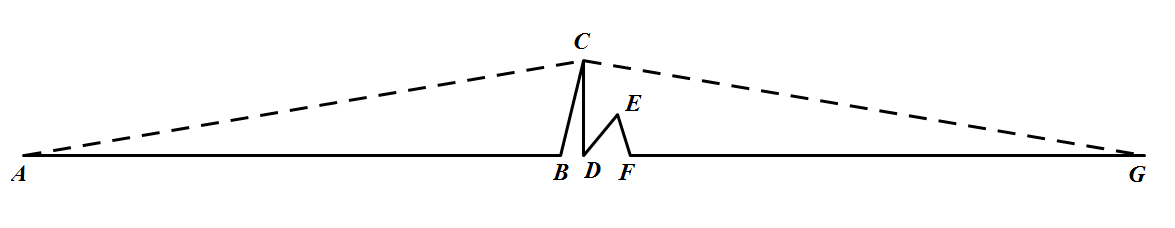} 
\caption{}\label{fig:4} 
\end{figure}

For a fixed $\tau$ with $1<\tau<1.001$, the corresponding basic figure is as in Figure~\ref{fig:4}. The points $B, D, F$ lie on segment $\ol{AG}$, and the magnitudes of the segments are $\ol{AG}=1$, $\ol{AD}=1/2$, $\ol{AB}=(7/15)^{1/\tau}$, $\ol{FG}=7/15$. The magnitudes of the angles are $\angle CAG=\angle CGA=\pi/18$, $\angle CDE=\pi/9$. The position of point $E$ is determined by $\ol{DE}=(7/15)^{1/\nu}\ol{CD}$, where
\baa{\label{52}\nu=
\begin{cases}
\tau, &\text{if $\tau$ is irrational,}\cr 
\tau-\frac{\tau-1}{\sqrt{2}}, &\text{if $\tau$ is rational.}\cr
\end{cases}}
Note that $\nu$ is always irrational and $1<\nu<1.001$.

Let $E'$ be the projection of $E$ on $\ol{AG}$. Then
\ba{\ol{DE'}&=(1/2)\tan(\pi/18)(7/15)^{1/\tau}\sin(\pi/9)\\
&<(1/2)(7/15)^{1000/1001}(1-\cos(\pi/9))\\
&<1/30.}
Thus $E'$ is between $D$ and $F$. It follows that $\angle EFG>\pi/2$ and $E$ is in the interior of triangle $CAG$. Therefore, polygon $ABCDEFG$ is a basic figure with basic triangle $\bigtriangleup CAG$. 

We denote polygon $ABCDEFG$ by $\Om_\tau$, and the corresponding  self-similar set by $\Lm_\tau$. The corner angles satisfy $\theta_j>0$ for $j\not=3$ and $\theta_3=0$. It is clear that $\eta_0>0$. 
By Proposition~\ref{proposition4.2}, in order to show that $\Lm_\tau$ is a regular self-similar arc, it suffices to verify
\ba{u\not\in \Sigma:=\{jx-ky: j,k\in\NN\},}
where
\baa{\label{123} x=\zeta,\;\; y=\tau^{-1}\zeta, \;\;\zeta=\log(15/7),\;\; u=\log\frac{\ol{CD}\cdot\ol{CG}}{\ol{DE}\cdot\ol{AG}}=\nu^{-1}\zeta.}
If $\tau$ is rational, then $\nu=\tau-(\tau-1)/\sqrt{2}$ is irrational, and for $j,k=0,1,2,\dots$,
\ba{u-(jx-ky)=\zeta(\nu^{-1}-(j-k/\tau))\not=0.}
If $\tau$ is irrational, then $\nu=\tau$, and for $j,k=0,1,2,\dots$,
\ba{u-(jx-ky)=\tau^{-1}\zeta j(-\tau+(k+1)/j)\not=0.}
Therefore $u\not\in\Sigma$ and $\Lm_\tau$ is a regular self-similar arc.

When $\tau$ is rational, $x/y$ is rational, hence $\Lm_\tau$ is a 1-quasi-arc.
We now consider the case where $\tau$ is irrational. In this case we have $\nu=\tau$ and $u=y$.
By Proposition~\ref{pro8.1}, $\Lm_\tau$ is a $t$-quasi-arc \iiff  there is a constant $M_{\tau,t}>0$ such that
\baa{\label{exp}e^{-j(t-1)x}\le M_{\tau,t}|u-jx+ky|}
In light of (\ref{123}), inequality (\ref{exp}) is reduced to
\ba{\tau(M_{\tau,t}\zeta)^{-1}j^{-1}e^{-j(t-1)\zeta}<  |\tau-(k+1)/j|, }
which is equivalent to 
\ba{C_{\tau,t}j^{-1}e^{-j(t-1)\zeta}<  |\tau-k/j|.}
Therefore, for  $t>1$, $\Lm_\tau$ is a $t$-quasi-arc \iiff $\tau$ has approximation property $J_{(t-1)\zeta}$.

We summarize the above discussion as follows.

\bE {\rm For $1<\tau<1.001$, let $\Lm_\tau$ be the regular self-similar arc generated by the basic figure $\Om_\tau$ in Figure~\ref{fig:4}, where
\ba{ &  \ol{AG}=1,\; \ol{AD}=1/2,\; \ol{AB}=(7/15)^{1/\tau},\; \ol{FG}=7/15,\\ &\angle CAG=\angle CGA=\pi/18,\; \angle CDE=\pi/9,\; \ol{DE}=(7/15)^{1/\nu}\ol{CD},}
where $\nu$ is defined by (\ref{52}). For  $t>1$, $\Lm_\tau$ is a $t$-quasi-arc \iiff $\tau$ has approximation property $J_{(t-1)\zeta}$.

(1) For a fixed $t_0>1$, by Theorem~\ref{7.11}, there is a transcendental number $\tau\in (1, 1.001)$ such that $\tau$ has approximation property $J_{(t_0-1)\zeta}$, but $\tau$ has approximation property  $J_{(t-1)\zeta}$ for no $t\in(1,t_0)$. For such a $\tau$, $\Lm_\tau$ is a $t_0$-quasi-arc, but $\Lm_\tau$ is a $t$-quasi-arc for no $t<t_0$. 

(2) For a fixed $t_0>1$, by Theorem~\ref{7.12}, there is a transcendental number $\tau\in (1, 1.001)$ such that  $\tau$ does not have approximation property $J_{(t_0-1)\zeta}$, but $\tau$ has approximation property $J_{(t-1)\zeta}$ for each $t>t_0$. For such a $\tau$, $\Lm_\tau$ is a $t$-quasi-arc for each $t>t_0$, but $\Lm_\tau$ is not a $t_0$-quasi-arc.

(3) By Theorem~\ref{7.13}, there is a transcendental number $\tau\in (1, 1.001)$ such that $\tau$ has approximation property $J_{(t-1)\zeta}$ for each $t>1$. For such a $\tau$, $\Lm_\tau$ is a $t$-quasi-arc for each $t>1$. Since $x/y=\tau$, which is an irrational number, it follows from   Theorem~\ref{thmrat} that $\Lm_\tau$ is not a $1$-quasi-arc.  

(4) By Theorem~\ref{7.14}, there is a transcendental number $\tau\in (1, 1.001)$ such that $\tau$ has approximation property $J_{(t-1)\zeta}$ for no $t>1$. For such a $\tau$, $\Lm_\tau$ is a $t$-quasi-arc for no $t> 1$. }\eE

As a consequence of the example, we obtain the following theorem.

\bT\label{7.15}
(1) There exists a regular self-similar arc $\Lm$ with $\theta_{\min}=0$ such that $\Lm$ is $t$-quasi-arc for no $t\ge 1$.

(2) Let $t_0> 1$. Then there exists a regular self-similar arc $\Lm$ with $\theta_{\min}=0$ such that $\Lm$ is a $t_0$-quasi-arc, but $\Lm$ is a $t$-quasi-arc for no $1\le t<t_0$.

(3) Let $t_0\ge 1$. Then there exists a regular self-similar arc $\Lm$ with $\theta_{\min}=0$ such that $\Lm$ is a $t$-quasi-arc for each $t>t_0$, but $\Lm$ is not a $t_0$-quasi-arc.
\eT

\brs When $t_0=1$, the third part of  Theorem~\ref{7.15} says that there exists a regular self-similar arc $\Lm$ with $\theta_{\min}=0$ which is an $t$-quasi-arc for each $t\in (1,s)$, where $s$ is the Hausdorff dimension of $\Lm$. Then by Theorem~ \ref{thm3}, the Hausdorff measure function $f$ is a Whitney function on $\Lm$.
\ers

{\bf Acknowledgment.} Part of the first named author's work was done while visiting Tshinghua University Yau Mathematical Sciences Center during his sabbatical leave
in spring 2014. He is grateful for the Center's hospitality and financial support.

%As consequence of Example~\ref{7.112}, we have the following corollary.
%\bC\label{10.2}There is a regular self-similar arc $\Lm$ with $\theta_{\min}=0$ such that $\Lm$ is a $t$-quasi-arc for each $t>1$, but $\Lm$ is not a 1-quasi-arc. 
%For such a $\Lm$, the Hausdorff measure function is a Whitney function.\eC

\end{document}